\documentclass[11pt]{article}

\usepackage[square,numbers]{natbib}
\usepackage{latexsym,mathrsfs}
\usepackage{amsmath,amssymb}
\usepackage{amsthm,enumerate,verbatim}
\usepackage{amsfonts}
\usepackage{graphicx}
\usepackage{algorithm}
\usepackage{algorithmicx}
\usepackage{algpseudocode}
\usepackage{url}
\usepackage[colorlinks=true,linkcolor=red, citecolor=green, urlcolor=blue]{hyperref}
\usepackage{booktabs}
\usepackage{subfig}
\usepackage{hyperref}
\usepackage{tikz}
\usepackage{authblk}
\usepackage[titletoc,title]{appendix}

\usepackage{lscape}
\usepackage{longtable}

\newtheorem{lemma}{Lemma}

\newtheorem{theorem}{Theorem}
\newtheorem{corollary}{Corollary}

\newtheorem*{definition*}{Definition}

\newtheorem{remark}{Remark}

\DeclareMathOperator{\argmin}{argmin}

\title{Solving the Maximum Clique Problem with Symmetric Rank-One Nonnegative Matrix Approximation}
\date{}

\author[1]{Melisew Tefera Belachew\thanks{This research was conducted when the first author was visiting the Department of Mathematics and Operational \hspace*{1.8em}Research at the University of Mons.}\thanks{Email: melisew.belachew@uniba.it;\ melisewt@gmail.com}}

\author[2]{Nicolas Gillis\thanks{Email: nicolas.gillis@umons.ac.be. NG acknowledges the support by the F.R.S.-FNRS, through the incentive grant for scientific research n$^\text{o}$ F.4501.16.}}

\affil[1]{Department of Mathematics, Universit\`{a} degli Studi di Bari Aldo Moro

Via E. Orabona~4, I-70125, Bari, Italy}

\affil[2]{Department of Mathematics and Operational Research,
Facult\'e Polytechnique,

Universit\'e de Mons, Rue de Houdain~9, 7000 Mons, Belgium}

\begin{document}

\maketitle

\begin{abstract}
Finding complete subgraphs in a graph, that is, cliques, is a key problem and has many real-world applications, e.g., finding communities in social networks, clustering gene expression data, modeling ecological niches in food webs, and describing chemicals in a substance.
The problem of finding the largest clique in a graph is a well-known $\mathcal{NP}$-hard problem and is called the maximum clique problem (MCP). In this paper, we formulate a very convenient continuous characterization of the MCP based on the symmetric rank-one nonnegative approximation of a given matrix,
and build a one-to-one correspondence between stationary points of our formulation and cliques of a given graph. In particular, we show that the local (resp.\@ global) minima of the continuous problem corresponds to the maximal (resp.\@ maximum) cliques of the given graph.
We also propose a new and efficient clique finding algorithm based on our continuous formulation and test it on various synthetic and real data sets to show that the new algorithm outperforms other existing algorithms based on the Motzkin-Straus formulation, and can compete with a sophisticated combinatorial heuristic. 

\end{abstract}

\section{Introduction}

A graph $G$ is a pair of sets $(V,E)$, where $V$ is the set of vertices and $E\subseteq V\times V$ is the set of edges formed by pairs of vertices. In this paper we only consider undirected and simple graphs. An undirected graph is complete if and only if all vertices are connected to one another, that is, $E = V\times V$. A \emph{clique} in an undirected graph $G$ is a subset of its vertices such that the corresponding subgraph is complete.
A clique is said to be \emph{maximal} if it is not contained in any larger clique.
A \emph{maximum} clique is a clique with maximum number of vertices (or, equivalently, edges) and the maximum clique problem (MCP) is the problem of finding such a clique. The number of vertices in a maximum clique in a graph $G$ is called the \emph{clique number} of $G$ and is denoted by $\omega (G)$. The MCP is a well-known $\mathcal{NP}$-hard problem and the associated decision problem, that is, the task of checking whether there is a clique of a given size in a graph, is $\mathcal{NP}$-complete~\cite{Gareyetal}.

The MCP arises in many real-life applications. The word ``clique'' in its graph-theoretic usage was first introduced in~\cite{Luceetal} where experts used complete graphs to model groups of people who all know each other in social networks.
The MCP can also be used to model the problem of clustering gene expression data in bioinformatics~\cite{BenDoretal},
to model ecological niches in food webs~\cite{Sugihara}, to analyze telecommunication networks~\cite{Prihar},
and to describe chemicals in a substance that have a high degree of similarity with a target structure~\cite{Rhodesetal}.

Since the MCP is an $\mathcal{NP}$-hard problem, it is a challenging task to devise algorithms to identify large cliques in graphs.
Experts use various approaches to tackle this problem among which the Motzkin-Straus continuous formulation is a well-known and widely used tool for MCP. There are also a multitude of discrete approaches for the MCP, which is one of the most fundamental problem in graph theory.
Performing a literature review of this extremely rich literature is out of the scope of this paper, and we refer the readers to~\cite{Bomze2, WuHao, Grossoetal} for surveys on these methods.

\subsection{Motzkin-Straus Formulation}

Let $G=(V,E)$ be an undirected graph where $V = \left\{1, 2, ..., n\right\}$ is the vertex set,
and let $A = \left(a_{ij}\right)_{i,j=1}^{n}$ be the binary adjacency matrix of $G$, where $a_{ij}=1$ if and only if $(i,j)\in E$ and $a_{ij}=0$ otherwise. The \emph{Motzkin-Straus formulation} of the MCP is given by~\cite{MotzkinStraus}
\begin{align}
\label{eq:MSMCP}
\max_{\mathbf{u}\in \mathbb{R}^{n}_{+}} \mathbf{u}^{\top}A\mathbf{u}
\textrm{\quad such that\quad }
\sum^{n}_{i=1}u_{i} = 1. \tag{MS}
\end{align}
\begin{theorem}[Motzkin and Straus~\cite{MotzkinStraus}]
\label{theorem1}
The optimal value of \eqref{eq:MSMCP} is given by
$1 - \frac{1}{\omega(G)}$, where $\omega(G)$ is the clique number of $G$.
\end{theorem}
As mentioned above, many exact and heuristic algorithms have been proposed to solve the MCP; see for example \cite{Bomze, Dingetal, Gibbonsetal, Pardalosetal2, Pelillo}.
Some of these methods use the above Motzkin-Straus formulation~\eqref{eq:MSMCP}. Note that using a continuous formulation of a combinatorial problem is a standard approach in global optimization; see, e.g., \cite{dgp97, hpt00} where it is used for the satisfiability problem.

An important result is that if the (zero) diagonal entries of $A$ are replaced by $\frac{1}{2}$, then any local (resp.\@ global) maximum of \eqref{eq:MSMCP} corresponds to a maximal (resp.\@ maximum) clique~\cite{Bomze}.
In this paper, we propose a new continuous formulation for the MCP based on a symmetric rank-one nonnegative matrix approximation problem.
Below, we give an overview of our continuous formulation and its advantages compared to~\eqref{eq:MSMCP}.

\subsection{Contribution and Outline of the Paper}

\label{sec:contriboutline}

Given a graph $G$ and its (symmetric) adjacency matrix $A$, let us define the associated modified adjacency matrix of $G$ by $B = A+I_{n}$, where $I_{n}$ is the identity matrix of dimension $n$.
Moreover, for some parameter $d\geq 0$, we define a (symmetric) matrix $M_{d}=\left(m_{ij}\right)_{i,j=1}^{n}$ as follows:
\begin{align}
\label{eq:matrixM}
\begin{split}
m_{ij} = \left\{ \begin{array}{l} \ \ 1, \textrm{\ if \ } b_{ij}=1, \\
-d, \textrm{\ if \ } b_{ij}=0.
\end{array} \right.
\end{split}
\end{align}

We propose in this paper to study the following optimization problem
\begin{equation}
\label{eq:MCP00}
\max_{\mathbf{u}\in \mathbb{R}^{n}_{+}} \ \mathbf{u}^\top M_{d} \mathbf{u} \textrm{\quad such that\quad }
\sum_{i=1}^{n}u_{i}^{2} \leq 1.
\end{equation}
We will show that the optimal value of~\eqref{eq:MCP00} is $\omega(G)$ (Corollary~\ref{optval}).
Note the similarity between the Motzkin-Straus formulation~\eqref{eq:MSMCP} and~\eqref{eq:MCP00}. An important difference is that the feasible set of~\eqref{eq:MCP00} is smooth, hence it is easier to project onto it.
This is an advantage for example when designing nonlinear optimization methods, e.g., projected gradient methods; see section~\ref{sec:clqfinalg}.

In section~\ref{sec:section2}, we will show that \eqref{eq:MCP00} is equivalent to the following symmetric rank-one matrix approximation problem (see Theorem~\ref{theorem2new} for a rigorous characterization)
\begin{align}
\label{eq:MCP000}
\min_{\mathbf{u}\in \mathbb{R}^{n}_{+}} \left\|M_{d}-\mathbf{uu}^{\top}\right\|_{F}^{2}.
\end{align}
Unlike other formulations of the MCP, we will draw very precise relationships between stationary points of~\eqref{eq:MCP000}, or, equivalently, \eqref{eq:MCP00}, and cliques of  the graph $G$.

Our theoretical result is therefore more complete than for the Motzkin-Straus formulation as we can also associate to any stationary point of~\eqref{eq:MCP000} a clique of $G$.

The paper is organized as follows.
In section~\ref{sec:formulationMCP} we introduce our continuous formulation~\eqref{eq:MCP000} for the MCP.
In section~\ref{sec:contxicmcp}, we show that
\begin{itemize}
\item[$\bullet$] the two formulations \eqref{eq:MCP00} and \eqref{eq:MCP000} are equivalent (Theorem~\ref{theorem2new}),
\item[$\bullet$] the local (resp.\@ global) minima of \eqref{eq:MCP000} coincide with the maximal (resp.\@ maximum) cliques of a given graph $G$ (Theorems~\ref{theorem3} and~\ref{theorem4}), and
\item[$\bullet$] every stationary point of the continuous optimization problem~\eqref{eq:MCP000} coincide with a feasible solution of the MCP (that is, a clique); see Theorem~\ref{thmstfs}.
\end{itemize}
In section \ref{sec:clqfinalg}, we propose a new and efficient clique finding algorithm based on our continuous formulation and show that the limit points of this algorithm coincide with cliques of the graph $G$. In addition, we present various experimental results that show competitiveness of the new algorithm compared to other clique finding algorithms.

\section{Continuous Characterization of the MCP using Symmetric Nonnegative Matrix Approximation}
\label{sec:section2}

In this section we derive a new continuous formulation of the MCP using symmetric rank-one nonnegative matrix approximation.

\subsection{Symmetric Rank-One Nonnegative Matrix Approximation and the MCP}
\label{sec:formulationMCP}

The so-called (discrete) \emph{vertex formulation} of the MCP~\cite{Pardalosetal} is given by
\begin{equation} \label{eq:MCP}
\max_{\mathbf{u} \in \left\{0,1\right\}^n} \; \sum_{i=1}^{n} u_{i}
\quad
\textrm{ such that }
\quad
u_{i}+u_{j} \leq 1+a_{ij} \ \forall i \neq j.
\end{equation}
The $i$th vertex belongs to a feasible solution of~\eqref{eq:MCP} if and only if $u_{i}=1$, otherwise $u_{i}=0$.
The constraint $u_{i}+u_{j} \leq 1+a_{ij}$ ensures that if there is no edge between the vertices $i$ and $j$, that is, if $a_{ij}=0$, then either $u_{i}=0$ or $u_{j}=0$. Hence, there is a one-to-one correspondence between the feasible solutions of~\eqref{eq:MCP} and the cliques of $G$.

The objective function of~\eqref{eq:MCP} can be rewritten as follows.
Observe that maximizing $\sum_{i} u_{i}$ reduces to requiring as many ones as possible in the vector $\mathbf{u}$,
which will lead to having more ones in the matrix $\mathbf{uu}^{\top}$, hence maximizing $\sum_{i,j} u_{i}u_{j}$.
Let $B = A + I_n$ and $\mathbf{u}$ be a feasible solution of \eqref{eq:MCP}.
Since $B$ and $\mathbf{u}$ are binary and  $u_{i}u_{j}\leq b_{ij} \ \forall i,j$, we have that
$ 
\sum_{i,j=1}^{n} u_{i}u_{j}  =  \sum_{i,j=1}^{n} (u_{i}u_{j})^{2} = \sum_{i,j=1}^{n} b_{ij}(u_{i}u_{j})^{2} = \|B\|^{2}_{F} - \left\|B-\mathbf{uu}^\top\right\|^{2}_{F},
$ 
where $\|\cdot\|_{F}$ is the Frobenius norm. In fact,
$
\left\|B-\mathbf{uu}^{\top}\right\|^{2}_{F} = \|B\|^{2}_{F} - 2\langle B, \mathbf{uu}^{\top}\rangle + \langle\mathbf{uu}^{\top}, \mathbf{uu}^{\top}\rangle, 
$
where $\langle \cdot, \cdot \rangle$ is the Frobenius inner product.
Hence, \eqref{eq:MCP} is equivalent to
\begin{equation}
\label{eq:MC(G)}
\min_{\mathbf{u} \in \left\{0,1\right\}^n}
\;
\left\|B-\mathbf{uu}^{\top}\right\|^{2}_{F}
\quad
\textrm{ such that }
\quad
u_{i}+u_{j} \leq 1+b_{ij} \ \forall i\neq j, \tag{MC}
\end{equation}
where the objective function is equal to the number of vertices outside the clique, and its minimization is therefore equivalent to maximizing the vertices contained in the clique. Hence, \eqref{eq:MC(G)} approximates $B$ via a symmetric rank-one binary approximation.

\subsection{New Continuous Formulation of the MCP}
\label{sec:contxicmcp}

We start off by defining the following problem: given an $n$-by-$n$ matrix $M\in \mathbb{R}^{n\times n}$, find its best rank-one nonnegative matrix approximation, that is, solve
\begin{align}
\label{eq:R1NC}
\min_{\mathbf{u}\in \mathbb{R}^{n}_{+}} \left\|M-\mathbf{uu}^{\top}\right\|_{F}^{2}. \tag{R1NM}
\end{align}
Given a parameter $d\geq 0$, a graph $G$, and its adjacency matrix $A\in \left\{0,1\right\}^{n\times n}$, we define the associated modified adjacency matrix by $B = A+I_{n}$ as before and a matrix $M_{d} = (1 + d)B - d\mathbf{1}_{n\times n}$ as in \eqref{eq:matrixM}, where $\mathbf{1}_{n\times n}$ is the $n$-by-$n$ matrix of ones.
Then, using the matrix $M_{d}$, we can define the following instance of \eqref{eq:R1NC}:
\begin{align}
\label{eq:R1NdCG}
\min_{\mathbf{u}\in \mathbb{R}^{n}_{+}} F(\mathbf{u}) = \left\|M_{d}-\mathbf{uu}^{\top}\right\|_{F}^{2}. \tag{R1N$_{d}$M}
\end{align}

In this paper, we analyze \eqref{eq:R1NdCG} as a continuous formulation of the MCP. It shares some similarities with the Motzkin-Straus formulation~\eqref{eq:MSMCP}. In fact, we have the following result.
\begin{theorem}
\label{theorem2new}
Let $M_{d}$ be as defined in \eqref{eq:matrixM} and consider the following optimization problem
\begin{equation}
\label{eq:MCPEY}
 \max_{\mathbf{v} \in \mathbb{R}^{n}_{+}} \; \mathbf{v}^{\top} M_{d} \mathbf{v}
\quad
\textrm{ such that }
\quad
\| \mathbf{v} \|^{2}_{2} \leq 1.
\end{equation}
Then, there is a one-to-one correspondence between the nontrivial stationary points of \eqref{eq:R1NdCG} and the stationary points of \eqref{eq:MCPEY} with positive objective function value. Precisely,
$\mathbf{u} = \left(\mathbf{v}^{\top}M_{d}\mathbf{v}\right)^{1/2} \mathbf{v}$ is a nontrivial stationary point of \eqref{eq:R1NdCG} if and only if $\mathbf{v} = \frac{\mathbf{u}}{\|\mathbf{u}\|_{2}}$ is a stationary point of \eqref{eq:MCPEY} with a positive objective function value (that is, $\mathbf{v}^{\top} M_{d} \mathbf{v} > 0$ hence $\mathbf{v} \neq 0$ is a stationary point with objective function value larger than zero). Moreover, for these stationary points, we have
$ 
\|M_{d}-\mathbf{uu}^{\top}\|^{2}_{F} = \|M_{d}\|^{2}_{F} - \left(\mathbf{v}^{\top}M_{d}\mathbf{v}\right)^{2}.
$ 
\end{theorem}
{\it Proof}
First, we derive the optimality conditions and stationary points of \eqref{eq:R1NdCG} and \eqref{eq:MCPEY}. The first-order optimality conditions of \eqref{eq:R1NdCG} are
\begin{equation}
\label{eq:optimalityclique}
\mathbf{u}\geq 0, \quad \nabla_{\mathbf{u}} F(\mathbf{u}) = \mathbf{uu}^{\top}\mathbf{u}-M_{d}\mathbf{u}\geq 0 \quad \textrm{\ and\ } \quad \mathbf{u}\odot \nabla_{\mathbf{u}} F(\mathbf{u}) = 0,
\end{equation}
where $\odot$ is the Hadamard product, that is, a vector $\mathbf{u}$ is a stationary point of \eqref{eq:R1NdCG} if and only if it satisfies the optimality conditions \eqref{eq:optimalityclique}. The optimality conditions given in \eqref{eq:optimalityclique} can be written equivalently as
\begin{equation}
\label{eq:optimalityclique2}
\mathbf{u} = 0 \quad \textrm{or} \quad \mathbf{u} = \max \left(0, \frac{M_{d}\mathbf{u}}{\|\mathbf{u}\|_{2}^{2}}\right) = \frac{\left[M_{d}\mathbf{u}\right]_{+}}{\|\mathbf{u}\|_{2}^{2}}.
\end{equation}
In fact, if $u_i = 0$, then $(M_{d}\mathbf{u})_i \leq 0$ (the gradient is nonnegative), while if  $u_i > 0$, then $u_i = (M_{d}\mathbf{u})_i$ (the gradient is equal to zero).
Hence, the nontrivial stationary points of \eqref{eq:R1NdCG} satisfy
\begin{equation}
\label{eq:stpoint1}
\mathbf{u} = \frac{[M_{d}\mathbf{u}]_{+}}{\|\mathbf{u}\|_{2}^{2}}.
\end{equation}
Next consider the Lagrangian function of \eqref{eq:MCPEY}
$
\mathcal{L}(\mathbf{v};\boldsymbol{\lambda},\mu) = -\frac{1}{2}\mathbf{v}^{\top}M_{d}\mathbf{v} +
\frac{\mu}{2}(\|\mathbf{v}\|^{2}_{2}-1) - \mathbf{v}^{\top}\boldsymbol{\lambda},
$
where $\boldsymbol{\lambda}\in \mathbb{R}^{n}_{+}$ and $0 \leq \mu \in \mathbb{R}$.
Any first-order (nontrivial) stationary point of \eqref{eq:MCPEY} with positive objective function value satisfies
\begin{equation}
\label{eq:opticondv}
\mathbf{v}^{\top}M_{d}\mathbf{v} > 0,
\ \nabla_{\mathbf{v}} \mathcal{L}= -M_{d}\mathbf{v} + \mu \mathbf{v} - \boldsymbol{\lambda} = 0,
\ \mu \geq 0,
\ \|\mathbf{v}\|^{2}_{2} \leq 1,
\ \mu (1- \|\mathbf{v}\|^{2}_{2} ) = 0,
\ \boldsymbol{\lambda}\geq0,
\ \mathbf{v}\geq 0,
\ \lambda_{i}v_{i} = 0, \ \forall i.
\end{equation}
If $\mu=0$, we have $M_{d}\mathbf{v} = -\boldsymbol{\lambda}\leq0$ hence  $\mathbf{v}^{\top}M_{d}\mathbf{v} = -\mathbf{v}^{\top} \boldsymbol{\lambda} \leq 0$, a contradiction.  Therefore, $\mu > 0$ so that the condition  $\mathbf{v}^{\top} M_{d}\mathbf{v} = \mu \mathbf{v}^{\top} \mathbf{v} = \mu > 0$ is satisfied, since $\|\mathbf{v}\|_2 = 1$ by complementarity.
Hence, similarly as for $\mathbf{u}$, the stationary points of \eqref{eq:MCPEY}  with positive objective function value are given by
\begin{equation}
\label{eq:stpoint2}
\mathbf{v} = \frac{M_{d}\mathbf{v} + \boldsymbol{\lambda}}{\mu} = \frac{[M_{d}\mathbf{v}]_{+}}{\mu}
=  \frac{[M_{d}\mathbf{v}]_{+}}{\| [M_{d}\mathbf{v}]_{+} \|_2}  ,
\end{equation}
since $\|\mathbf{v}\|_2 = 1$ (and $\mu > 0$). The second equality comes from the fact that $v_{i} \lambda_{i} = 0$ for all $i$.

We can now prove that the conditions \eqref{eq:stpoint1} and \eqref{eq:stpoint2} coincide when $\mathbf{u} = (\mathbf{v}^{\top}M_{d}\mathbf{v})^{1/2} \mathbf{v}$ and $\mathbf{v} = \frac{\mathbf{u}}{\|\mathbf{u}\|_{2}}$:
\begin{itemize}

\item[$\bullet$] Let $\mathbf{u}$ satisfy~\eqref{eq:stpoint1}, then $\mathbf{v} = \frac{\mathbf{u}}{\|\mathbf{u}\|_{2}}$
satisfies~\eqref{eq:stpoint2}: in fact,
\begin{equation}
\label{eq:stpoint3}
\mathbf{v} = \frac{\mathbf{u}}{\|\mathbf{u}\|_{2}} \stackrel{\eqref{eq:stpoint1}}{=} \frac{[M_{d}\mathbf{u}]_{+}}{\|\mathbf{u}\|_{2}^{2}} \frac{1}{\frac{\|[M_{d}\mathbf{u}]_{+}\|_{2}}{\|\mathbf{u}\|_{2}^{2}}}  = \frac{[M_{d}\frac{\mathbf{u}}{\|\mathbf{u}\|_{2}}]_{+}}{\|[M_{d}\frac{\mathbf{u}}{\|\mathbf{u}\|_{2}}]_{+}\|_{2}}
 = \frac{[M_{d}\mathbf{v}]_{+}}{\|[M_{d}\mathbf{v}]_{+}\|_{2}}.
\end{equation}
\item[$\bullet$] Let $\mathbf{v}$ satisfy~\eqref{eq:stpoint2}, then $\mathbf{u} = \left(\mathbf{v}^{\top}M_{d}\mathbf{v}\right)^{1/2}\mathbf{v}$ satisfies~\eqref{eq:stpoint1}. In fact, since $\mu = \mathbf{v}^{\top}M_{d}\mathbf{v}$ and $\|\mathbf{v}\|^{2}_{2}=1$, we have
\begin{equation}
\label{eq:stpointu1}
\mathbf{u}
= \left(\mathbf{v}^{\top}M_{d}\mathbf{v}\right)^{1/2}\mathbf{v} \stackrel{\eqref{eq:stpoint2}}{=} \left(\mathbf{v}^{\top}M_{d}\mathbf{v}\right)^{1/2}\frac{[M_{d}\mathbf{v}]_{+}}{\mu}
= \frac{[M_{d}\left(\mathbf{v}^{\top}M_{d}\mathbf{v}\right)^{1/2}\mathbf{v}]_{+}}{\mu}
= \frac{[M_{d}\mathbf{u}]_{+}}{\mathbf{v}^{\top}M_{d}\mathbf{v}}
= \frac{[M_{d}\mathbf{u}]_{+}}{\|\mathbf{u}\|^{2}_{2}},
\end{equation}
where the last equality is due to $\|\mathbf{u}\|^{2}_{2} = \|\left(\mathbf{v}^{\top}M_{d}\mathbf{v}\right)^{1/2}\mathbf{v}\|^{2}_{2} = \mathbf{v}^{\top}M_{d}\mathbf{v}$.

\end{itemize}

Note that the scaling factor $\mu = \mathbf{v}^{\top}M_{d}\mathbf{v}$ is such that
$\mathbf{u} \mathbf{u}^{\top} = \left(\mathbf{v}^{\top}M_{d}\mathbf{v}\right) \mathbf{v} \mathbf{v}^{\top}$ approximates $M_{d}$ as well as possible: in fact,
$
\argmin_{\mu \geq 0} \|M_{d}- \mu \mathbf{vv}^{\top}\|^{2}_{F} = \max\left(0, \mathbf{v}^{\top}M_{d}\mathbf{v}\right).
$
Finally, for these nontrivial stationary points, we have that
\[
\|M_{d} - \mathbf{uu}^{\top}\|^{2}_{F}
= \|M_{d}\|^{2}_{F} - 2\mathbf{u}^{\top}M_{d}\mathbf{u} + \|\mathbf{u}\|^{4}_{2}  \\
= \|M_{d}\|^{2}_{F} - 2 (\mathbf{v}^{\top}M_{d}\mathbf{v})^2 + (\mathbf{v}^{\top}M_{d}\mathbf{v})^2 \|\mathbf{v}\|^{4}_{2}  \\
= \|M_{d}\|^{2}_{F} - \left(\mathbf{v}^{\top}M_{d}\mathbf{v}\right)^{2},
\]
which concludes the proof.
\qed

\begin{remark}
We have to consider nontrivial stationary point of~\eqref{eq:R1NdCG} because $\mathbf{u}=0$ is always stationary for~\eqref{eq:R1NdCG} while it does not correspond to a feasible solution of \eqref{eq:MCPEY} with unit norm.
The reason no to consider stationary points $\mathbf{v}$ of~\eqref{eq:MCPEY} with negative objective function value is that
$\mathbf{u} = \left(\mathbf{v}^{\top}M_{d}\mathbf{v}\right)^{1/2} \mathbf{v}$  can  locally be improved simply be taking it closer to zero (e.g., multiply it by any $0 < \alpha < 1$ as $\mathbf{u}$ is not required to have unit norm) hence it is not a stationary point of~\eqref{eq:R1NdCG}.
For example, the rank-one matrix $D=-zz^T$ for some vector $z$ with unit norm has one negative eigenvalue (-1), and the optimal solution of~\eqref{eq:R1NC} is $\mathbf{u}=0$ (note that this is true even if $\mathbf{u}$ is not required to be nonnegative).
\end{remark}

\paragraph{Interpretation of~\eqref{eq:R1NdCG} and organization of the section.}
The parameter $d$ in~\eqref{eq:R1NdCG} can be interpreted as a penalty parameter that is used to satisfy the constraint
$u_i u_j \leq b_{ij}$ for all $i,j$.
In fact, the $-d$ entries in the matrix $M_{d}$ penalize the fact that entries in $\mathbf{uu}^{\top}$ are positive when corresponding to the zero entries of $B$: for each $i \neq j$ such that $b_{ij} = 0$, the term in the objective function is
$
(-d - u_iu_j)^2 = d^2 + 2 d u_i u_j + (u_i u_j)^2.
$
Therefore, as $d$ increases, the nondiagonal entries of $\mathbf{uu}^{\top}$ corresponding to the zero entries of $B$ are encouraged to be closer to zero (see Lemma~\ref{lemmad}).

In the next sections, we show that for $d \geq n$, the local (resp.\@ global) minima of the continuous optimization problem~\eqref{eq:R1NdCG} are binary and coincide with the maximal (resp.\@ maximum) cliques of the graph $G$, respectively (Theorems~\ref{theorem3} and~\ref{theorem4});
or equivalently, with the optimal solutions of the discrete problem \eqref{eq:MC(G)}.
Moreover, we show that the other stationary points of~\eqref{eq:R1NdCG} get arbitrarily close to the cliques of $G$ as $d$ increases
(Theorem~\ref{thmstfs}).
First, we explain the connections between the results of this paper, and the results from~\cite{GG14} on the the maximum-edge biclique problem.

\subsubsection{Link with the Maximum-Edge Biclique Problem and the Results from~\cite{GG14}}

Given a bipartite graph $\hat{G}$, the maximum-edge biclique problem (MBP) is the problem of finding a complete subgraph (that is, a biclique) with maximum number of edges.
Precisely, let $\hat{A}\in \{0,1\}^{m\times n}$ be the binary biadjacency matrix of $\hat{G}$. Similarly as the clique problem~\eqref{eq:MCP}, the MBP can be formulated as follows
\[
\min_{\mathbf{u}\in \{0,1\}^{m}, \mathbf{v}\in \{0,1\}^{n}} \;
\|\hat{A} - \mathbf{u}\mathbf{v}^{\top}\|^{2}_{F}
 \quad
\text{ such that }
\quad
u_{i} + v_{j} \leq 1 + \hat{a}_{ij} \ \forall i,j.
\]
Defining $\hat{M}_{d} = (1+\hat{d})\hat{A} - \hat{d}\mathbf{1}_{m\times n}$, where $\hat{d}$ is a positive parameter and $\mathbf{1}_{m\times n}$ is an $m$-by-$n$ matrix of ones, Gillis and Glineur~\cite{GG14} proposed the following continuous formulation of the MBP, a rank-one matrix approximation problem,
\begin{equation}
\label{eq:R1NdGB0}
\min_{\mathbf{u}\in \mathbb{R}^{m}_{+},\mathbf{v}\in \mathbb{R}^{n}_{+}}\|\hat{M}_{d}-\mathbf{u}\mathbf{v}^{\top}\|^{2}_{F}.
\end{equation}
They proved a one-to-one correspondence between the maximal (resp.\@ maximum) bicliques of the graph $\hat{G}$ and the local (resp.\@ global) minima of~\eqref{eq:R1NdGB0} for $\hat{d} \geq \max(m,n)$. They also showed that, as $d$ increases, all stationary points of~\eqref{eq:R1NdGB0} get arbitrarily close to the bicliques of $\hat{G}$.

Interestingly, if $\mathbf{u}$ is a stationary point of~\eqref{eq:R1NdCG}, then $(\mathbf{u},\mathbf{u})$ is a stationary point of~\eqref{eq:R1NdGB0} for the same matrix $\hat{M}_{d} = M_{d}$ (Lemma~\ref{lemma4}),
therefore we will be able to use the result from~\cite{GG14} to prove that
$\mathbf{u}$ gets closer to a clique of $G$ as $d$ increases (Theorem~\ref{thmstfs}).
However, the link between the maximal (resp.\@ maximum) cliques of $G$ and the local (resp.\@ global) minima of~\eqref{eq:R1NdCG} does not follow directly from~\cite{GG14}. For example, for the modified adjacency matrix
\[
B = \left( \begin{array}{ccccc}
1 & 1 & 1 & 1 & 1 \\
1 & 1 & 0 & 0 & 0 \\
1 & 0 & 1 & 0 & 0 \\
1 & 0 & 0 & 1 & 0 \\
1 & 0 & 0 & 0 & 1 \\
\end{array} \right) ,
\]
the maximum cliques contain any two vertices, while the maximum bicliques are $(1,0,0,0,0) \times (1,1,1,1,1)$ and $(1,1,1,1,1) \times (1,0,0,0,0)$.
Despite these differences, our proofs will use some arguments from~\cite{GG14},
and we will follow a similar organization to prove the one-to-one correspondence between the maximal (resp.\@ maximum) cliques of the graph ${G}$ and the local (resp.\@ global) minima of~\eqref{eq:R1NdCG}.

\subsubsection{Definitions and Notations}

Let us introduce the definitions and notations that will be used to prove the main results of this paper.

A ball centered at $\mathbf{x}\in \mathbb{R}^{n}_{+}$ with radius $r$ and intersected with the nonnegative orthant is defined as
\begin{equation*}
\mathcal{B}_{+}(\mathbf{x},r) = \{\mathbf{y}\in \mathbb{R}_{+}^{n}\mid \|\mathbf{x}-\mathbf{y}\|_{2}\leq r\}.
\end{equation*}
A vector $\mathbf{u}$ is a local minimum of \eqref{eq:R1NdCG} if and only if there exists an $\epsilon > 0$ such that for all $\mathbf{v}\in \mathcal{B}_{+}(\mathbf{u},\epsilon)$, we have $\|M_{d}-\mathbf{uu}^{\top}\|_{F}^{2} \leq \|M_{d}-\mathbf{vv}^{\top}\|_{F}^{2}$. A vector $\mathbf{u}$ is a global minimum of \eqref{eq:R1NdCG} if and only if $\|M_{d}-\mathbf{uu}^{\top}\|_{F}^{2} \leq \|M_{d}-\mathbf{vv}^{\top}\|_{F}^{2}$ for all $\mathbf{v}\in \mathbb{R}^{n}_{+}$.

Given a positive real number $d$, we define the following three sets of vectors:
\begin{itemize}
\item[$\bullet$] $\textsf{S}_{\textsf{p}}$, corresponding to the set of nontrivial stationary points of \eqref{eq:R1NdCG}, that is,
\begin{equation*}
\textsf{S}_{\textsf{p}} = \{\mathbf{u}\in \mathbb{R}^{n}_{+}\mid \mathbf{u} \textrm{\ satisfies\ } \eqref{eq:optimalityclique2} \textrm{\ and\ } \mathbf{u}\neq 0 \}.
\end{equation*}
\item[$\bullet$] $\textsf{L}_{\textsf{m}}$, corresponding to the set of nontrivial local minima of \eqref{eq:R1NdCG}.
\item[$\bullet$] $\textsf{G}_{\textsf{m}}$, corresponding to the set of nontrivial global minima of \eqref{eq:R1NdCG}.
\end{itemize}
By definition, $\textsf{G}_{\textsf{m}} \subseteq \textsf{L}_{\textsf{m}} \subseteq \textsf{S}_{\textsf{p}}$.

Let us also define the following three sets of binary vectors:
\begin{itemize}
\item[$\bullet$] $\textsf{F}_{\textsf{s}}$, corresponding to the set of feasible solutions of \eqref{eq:MC(G)}, that is,
\begin{equation*}
\textsf{F}_{\textsf{s}} = \{\mathbf{u}\in \mathbb{R}^{n}_{+}\mid \mathbf{u} \textrm{\ is a feasible soultion of\ }\eqref{eq:MC(G)}\}.
\end{equation*}
\item[$\bullet$] \textsf{C}$_{\textsf{m}}$, corresponding to the maximal cliques of $G$, that is, $\mathbf{u}\in \textsf{C}_{\textsf{m}}$ if and only if $\mathbf{u}\in \textsf{F}_{\textsf{s}}$ and  $\mathbf{u}$ corresponds to a maximal clique of $G$.
\item[$\bullet$] \textsf{C}$_{\textsf{M}}$, corresponding to the maximum cliques of $G$, that is, $\mathbf{u}\in \textsf{C}_{\textsf{M}}$ if and only if $\mathbf{u}\in \textsf{F}_{\textsf{s}}$ and  $\mathbf{u}$ corresponds to a maximum clique of $G$.
\end{itemize}
By definition, $\textsf{C}_{\textsf{M}} \subseteq \textsf{C}_{\textsf{m}} \subseteq \textsf{F}_{\textsf{s}}$.

\subsubsection{Key Lemmas}

Given an $n$-by-$n$ symmetric matrix $M\in \mathbb{R}^{n\times n}$, its best symmetric rank-one approximation can be obtained by solving the following unconstrained minimization problem
\begin{equation}
\label{eq:R1UMC}
\min_{\mathbf{u}\in \mathbb{R}^{n}} \left\|M-\mathbf{uu}^{\top}\right\|^{2}_{F}. \tag{R1U}
\end{equation}
In this section we will prove various results regarding~\eqref{eq:R1UMC} that will be crucial for the remainder of the paper. The following lemma is well-known although we do not know the reference for the original proof; see, e.g.,~\cite[Th.~1.14]{diep}. We give the proof here for completeness.
\begin{lemma}
\label{lemma1}
The local minima of \eqref{eq:R1UMC} are global minima. All other nontrivial stationary points are either saddle points or local maxima.
\end{lemma}
{\it Proof}
Following exactly the same argument as in Theorem~\ref{theorem2new}, we can show that the stationary points of \eqref{eq:R1UMC} correspond to the stationary points of
\begin{equation}
\label{eq:utMu}
\max_{\mathbf{v}\in \mathbb{R}^{n}} \mathbf{v}^{\top}M\mathbf{v} \textrm{\quad such that\quad } \|\mathbf{v}\|^{2}_{2} \leq 1,
\end{equation}
with a positive objective function value. The only difference with the proof of Theorem~\ref{theorem2new} is that the nonnegativity constraints must be discarded (hence removing $[.]_+$ and the Lagrangian multipliers $\boldsymbol{\lambda}$ from the proof).

It is well-known that the stationary points of~\eqref{eq:utMu} with a nonnegative objective function value are given by the normalized eigenvectors of $M$ associated with the nonnegative eigenvalues. In particular, the optimal (nontrivial) solution of \eqref{eq:utMu} is the normalized eigenvector of $M$ associated with the largest positive eigenvalue.
Let us denote
$\lambda_1(M) \geq \lambda_2(M)  \geq \dots \geq \lambda_n(M)$ the eigenvalues of $M$ in nonincreasing order.
Let $\mathbf{v}_{1}$ be an optimal solution of \eqref{eq:utMu} with $\mathbf{v}_{1}^{\top}M\mathbf{v}_{1}=\lambda_{1} > 0$ and $\|\mathbf{v}_1\|_2=1$.
In addition,
let $\mathbf{x}$ be an arbitrary stationary point of \eqref{eq:R1UMC} such that $\mathbf{x}^{\top}M\mathbf{x}=\lambda < \lambda_{1}$.
Since $\mathbf{v}_{1}$ and $\mathbf{x}$ are associated with different eigenvalues of $M$, they are orthogonal.
Now define $\mathbf{w} = \sqrt{1-\epsilon^{2}}\mathbf{x} + \epsilon \mathbf{v}_{1}$ for $-1\leq \epsilon \leq 1$.
We have
\[
\|\mathbf{w}\|_2^2
= (\sqrt{1-\epsilon^{2}}\mathbf{x} + \epsilon \mathbf{v}_{1})^{\top} (\sqrt{1-\epsilon^{2}}\mathbf{x} + \epsilon \mathbf{v}_{1})
= (1-\epsilon^{2}) \|\mathbf{x}\|_2 + \epsilon^2 \|\mathbf{v}_1\|_2 \leq 1, \; \; \text{ and }
\]
\begin{align*}
\mathbf{w}^{\top}M\mathbf{w} &
=  (\sqrt{1-\epsilon^{2}}\mathbf{x} + \epsilon \mathbf{v}_{1})^{\top} M (\sqrt{1-\epsilon^{2}}\mathbf{x} + \epsilon \mathbf{v}_{1})
= (1-\epsilon^{2})\mathbf{x}^{\top}M\mathbf{x} + \epsilon^{2}\mathbf{v}_{1}^{\top}M\mathbf{v}_{1} \\
& =  \mathbf{x}^{\top}M\mathbf{x} + \epsilon^{2}(\mathbf{v}_{1}^{\top}M\mathbf{v}_{1}-\mathbf{x}^{\top}M\mathbf{x})
= \lambda + \epsilon^{2}(\lambda_{1}-\lambda)
> \lambda
= \mathbf{x}^{\top}M\mathbf{x}.
\end{align*}
Hence, $\mathbf{x}$ is not a local minimum. Therefore, it is either a saddle point or a local maximum.
\qed

\begin{lemma}
\label{lemma2}
For the (symmetric) matrix $M_{d}$ defined in~\eqref{eq:matrixM} with at least one entry equal to $-d$ with $d\geq n$,
any optimal solution $\mathbf{u}$ of~\eqref{eq:R1UMC} with $M = M_{d}$ is such that $\mathbf{u}$ contains at least one nonpositive entry.
\end{lemma}
{\it Proof}
Suppose $\mathbf{u}$ is an optimal solution of \eqref{eq:R1UMC} such that $\mathbf{uu}^{\top} > 0$.
Note that the diagonal entries of $M_{d}$ are equal to one and that $M_{d}$ contains at least one $-d$ entry, hence $M_{d}$ contains at least two $-d$ entries.
Let $r$ denote the number of entries equal to $-d$ in $M_{d}$ with $r\geq 2$.
Therefore, since $\mathbf{u} > 0$,
\[
\|M_{d}-\mathbf{uu}^{\top}\|^{2}_{F} > rd^{2}.
\]
By assumption, there exists $(i,j)$ such that $i \neq j$ and $m_{ij}=-d$, and, by symmetry, $m_{ji}=-d$.
Consider the vector $\mathbf{v}\in \mathbb{R}^{n}$
such that
$v_{i} =\sqrt{\frac{d}{2}}, v_{j} = -\sqrt{\frac{d}{2}}$
and
$v_{k} = 0 \ \forall k \neq i, j$.
We have that $v_iv_i = v_jv_j = \frac{d}{2}, v_jv_i = v_iv_j = -\frac{d}{2}$
and $v_k v_l = 0 \ \forall k,l \neq i,j$, hence
\begin{align*}
\|M_{d}-\mathbf{vv}^{\top}\|^{2}_{F} & = 2\left(\frac{d}{2} - 1 \right)^{2} + 2\left(\frac{d}{2}\right)^{2} + (r-2)d^{2} + n^{2} - (r+2) \\
& < 
 rd^{2} + n^{2} - d^{2} - (r+2) \stackrel{r \geq 2}{\leq} rd^{2} + n^{2} - d^{2} \stackrel{d \geq n}{\leq} rd^{2},
\end{align*}
a contradiction. Therefore, any optimal solution $\mathbf{u}$ of \eqref{eq:R1UMC} must contain at least one nonpositive entry.
\qed

\subsubsection{Local and Global Minima of \eqref{eq:R1NdCG}}

In this section we characterize the relationship between the local (resp.\@ global) minima of \eqref{eq:R1NdCG} and the maximal (resp.\@ maximum) cliques of a given graph $G$.

\begin{lemma}
\label{lemma3}
Let $G=(V,E)$ be a graph with at least one edge,
with binary adjacency matrix $A$ and modified adjacency matrix $B=A+I_{n}$,
and $d\geq n$.
Then $\textsf{L}_{\textsf{m}}\subseteq \textsf{C}_{\textsf{m}}$.
\end{lemma}
{\it Proof}
Let $M_{d}\in \{-d,1\}^{n\times n}$ be the matrix defined in \eqref{eq:matrixM}.  Let $\mathbf{u}\in \textsf{L}_{\textsf{m}}$. To show $\textsf{L}_{\textsf{m}}\subseteq \textsf{C}_{\textsf{m}}$, we need to show that $\mathbf{u}$ is a feasible solution of \eqref{eq:MC(G)} and $\mathbf{u}$ corresponds to a maximal clique of $G$.
The support of a vector $\mathbf{u}$ is defined as the set of indices corresponding to the nonzero entries of $\mathbf{u}$. Let us denote the (non-empty) index set of the support of $\mathbf{u}$ by $S$ and define $\mathbf{u}'=\mathbf{u}(S)$ and $M_{d}'=M_{d}(S,S)$ to be the subvector and the submatrix with indices in $S$ and $S\times S$, respectively. Let us also define $G'$ as the graph whose modified adjacency matrix is given by $B'=B(S,S)$. Since $\mathbf{u}$ is a local minimum of \eqref{eq:R1NdCG} and the objective functions of \eqref{eq:R1NdCG} and \eqref{eq:R1NCG'} differ only by a constant (and also $\mathbf{u}^{\top}M_{d}\mathbf{u}=\mathbf{u}'^{\top}M_{d}'\mathbf{u}'$), we have that $\mathbf{u}'$ is a local minimum of \eqref{eq:R1NCG'}
\begin{align}
\label{eq:R1NCG'}
\min_{\mathbf{u}'\in  \mathbb{R}^{|S|}_+} \left\|M_{d}'-\mathbf{u}'\mathbf{u}'^{\top}\right\|^{2}_{F}. \tag{R1NCG$'$}
\end{align}
To show that $\mathbf{u}$ is a feasible solution of~\eqref{eq:MC(G)}, we first suppose there is a $-d$ entry in $M_{d}'$.
Since $\mathbf{u}'$ is positive, it is located in the interior of the feasible domain of \eqref{eq:R1NCG'}. Therefore, it is a local minimum of the unconstrained problem \eqref{eq:R1UCM'}
\begin{align}
\label{eq:R1UCM'}
\min_{\mathbf{u}'\in \mathbb{R}^{|S|}}\left\|M_{d}'-\mathbf{u}'\mathbf{u}'^{\top}\right\|^{2}_{F}. \tag{R1UCM$'$}
\end{align}
Thus, Lemma~\ref{lemma1} implies that $\mathbf{u}'$ is a global minimum of \eqref{eq:R1UCM'}. Moreover, since $M_{d}'$ contains at least one $-d$ entry, Lemma~\ref{lemma2} asserts that $\mathbf{u}'$ contains a non-positive entry, a contradiction. Therefore, $M_{d}'$ does not contain a $-d$ entry and as a result we have $M_{d}' = \mathbf{1}_{|S|\times |S|}$. Since $\mathbf{u}'$ is a global minimum of \eqref{eq:R1UCM'} and $M_{d}' = \mathbf{1}_{|S|\times |S|}$, we must have $\mathbf{u}'\mathbf{u}'^{\top} = M_{d}' = \mathbf{1}_{|S|\times |S|}, \ \mathbf{u}' = \mathbf{1}_{|S|}$ and $\mathbf{u}$ is binary. Therefore, $\mathbf{u}$ is a feasible solution of \eqref{eq:MC(G)}, that is, $\mathbf{u}\in \textsf{F}_{\textsf{s}}$.

Finally, let us show that $\mathbf{u}$ corresponds to a maximal clique of $G$. Assume that $\mathbf{u}$ corresponds to a clique of $G$ which is not maximal, that is, assume without loss of generality that $\exists i \notin S$ such that $\mathbf{u}+\mathbf{e}_{i}$ corresponds to a larger clique of $G$ where $\mathbf{e}_{i}$ is the unit vector whose $i$-th entry is equal to one.
For any $0<\epsilon \leq 1$, let $\mathbf{v} = \mathbf{u}+\epsilon\mathbf{e}_{i}$. Then we have $\left\|M_{d}-\mathbf{v}\mathbf{v}^{\top}\right\|^{2}_{F}\leq \left\|M_{d}-\mathbf{u}\mathbf{u}^{\top}\right\|^{2}_{F}$, since the entries of $M_{d}$ corresponding to edges contained only in the larger clique $\{i\}\times S$ are $1$'s and are approximated by values between $0$ and $1$ in $\mathbf{v}\mathbf{v}^{\top}$ whereas they are approximated by zeros in $\mathbf{u}$.
This contradicts the assumption that $\mathbf{u}$ is a local minimum.
Hence, $\mathbf{u}$ must correspond to a maximal clique of $G$, that is, $\mathbf{u}\in \textsf{C}_{\textsf{m}}$.
\qed

The next result shows that all the maximal cliques of a given graph $G$ correspond to the local minima of \eqref{eq:R1NdCG}.
\begin{theorem}
\label{theorem3}
If $G$ is a graph with at least one edge and $d\geq n$, then $\textsf{C}_{\textsf{m}}=\textsf{L}_{\textsf{m}}$.
\end{theorem}
{\it Proof}
Since $d\geq n$, by Lemma~\ref{lemma3} we have $\textsf{L}_{\textsf{m}} \subseteq \textsf{C}_{\textsf{m}}$. We show that $\textsf{C}_{\textsf{m}}\subseteq \textsf{Lm}$ in Appendix~\ref{proofth3}.
\qed

The next result states the strong relationship between the global solutions of~\eqref{eq:R1NdCG} and the maximum cliques of a graph $G$.
\begin{theorem}
\label{theorem4}
If $G$ is a graph with at least one edge and $d\geq n$, then $\textsf{G}_{\textsf{m}}=\textsf{C}_{\textsf{M}}$.
\end{theorem}
{\it Proof}
Let $\mathbf{u}\in \textsf{G}_{\textsf{m}}$. Then, by definition, $\mathbf{u}\in \textsf{L}_{\textsf{m}}$ and by Theorem~\ref{theorem3} $\mathbf{u}\in \textsf{C}_{\textsf{m}}$. Hence, $\mathbf{u}$ is binary. Next, observe that the objective functions of \eqref{eq:R1NdCG} and \eqref{eq:MC(G)} differ only by a constant:
\begin{equation*}
\left\|M_{d}-\mathbf{u}\mathbf{u}^{\top}\right\|^{2}_{F} = \left\|B - \mathbf{u}\mathbf{u}^{\top}\right\|^{2}_{F} + (n^{2}-\|B\|_{F}^{2})d^{2}.
\end{equation*}
Therefore, $\mathbf{u} \in \textsf{G}_{\textsf{m}}$ if and only if $\mathbf{u}\in \textsf{C}_{\textsf{M}}$.
\qed

\begin{corollary} \label{optval}
The optimal value of~\eqref{eq:R1NdCG} is $\|M_{d}\|_F^2 - \omega(G)^2$, and the optimal value of \eqref{eq:MCP00} is $\omega(G)$.
\end{corollary}
{\it Proof}
This follows from Theorem~\ref{theorem2new} (equivalence of~\eqref{eq:R1NdCG} and~\eqref{eq:MCP00}) and Theorem~\ref{theorem4} (an optimal solution of~\eqref{eq:R1NdCG} is the indicator vector corresponding to a maximum clique). In fact, we have the following: Let $\mathbf{u}$ and $\mathbf{v}$ be the global minima of \eqref{eq:R1NdCG} and \eqref{eq:MCP00}, respectively. Then,
\begin{align*}
\|M_{d}-\mathbf{uu}^{\top}\|^{2}_{F} \stackrel{\text{Thm}.~\ref{theorem2new}}{=} \|M_{d}\|^{2}_{F} - \left(\mathbf{v}^{\top}M_{d}\mathbf{v}\right)^{2} = \|M_{d}\|^{2}_{F} - \|\mathbf{u}\|^{4}_{2} = \|M_{d}\|^{2}_{F} - \omega (G)^{2},
\end{align*}
where the last equality is because Theorem~\ref{theorem4} implies that the global minima of \eqref{eq:R1NdCG} is binary, as a result we have $\|\mathbf{u}\|^{2}_{2} = \omega (G) = \mathbf{v}^{\top}M_{d}\mathbf{v}$.
\qed

\begin{corollary}
\eqref{eq:R1NdCG} is $\mathcal{NP}$-hard.
\end{corollary}
{\it Proof}
By Theorem~\ref{theorem4} finding the global optima of~\eqref{eq:R1NdCG} is equivalent to solving \eqref{eq:MC(G)} which is equivalent to solving \eqref{eq:MSMCP}, a well-known $\mathcal{NP}$-hard problem~\cite{Gareyetal}.
\qed

\subsubsection{Stationary Points and Maximal Cliques}

Here, we prove an important result that tells us how the maximal cliques of a given graph $G$ are related to the stationary points of \eqref{eq:R1NdCG} and to the feasible solutions of \eqref{eq:MC(G)}.
\begin{theorem}
\label{thoerem4}
If $G$ is a graph with at least one edge and $d\geq n$, then $\textsf{C}_{\textsf{m}}=\textsf{F}_{\textsf{s}}\cap\textsf{S}_{\textsf{p}}$.
\end{theorem}
{\it Proof}
Here we need to show that if $\mathbf{u}\in \textsf{C}_{\textsf{m}}$ then $\mathbf{u}$ belongs to both $\textsf{F}_{\textsf{s}}$ and $\textsf{S}_{\textsf{p}}$. We know that if $\mathbf{u}\in \textsf{C}_{\textsf{m}}$ then it is automatically binary. Let $S$ denote the non-empty support of $\mathbf{u}$. By definition, we have that $\mathbf{u}\in \textsf{F}_{\textsf{s}}$ and $\mathbf{u}$ corresponds to a maximal clique of $G$, that is,
\begin{equation}
\label{eq:uid}
\nexists i \textrm{\ such\ that\ }u_{i}=0 \textrm{\ and\ } m_{ij} = 1 \ \forall j \in S.
\end{equation}
It remains to show that $\mathbf{u}\in \textsf{S}_{\textsf{p}}$. For all $i$ such that $u_{i}=0$, by \eqref{eq:uid} at least one entry of $M_{d}(i,:)$ is $-d$. Therefore, we have
\begin{equation*}
\label{eq:uid2}
M_{d}(i,:)\mathbf{u}\leq(n-2)-d \stackrel{d\geq n}{<} 0.
\end{equation*}
Since $\mathbf{u}$ is binary we also have
\begin{equation}
\label{eq:uiM}
u_{i} = 0 \ \textrm{and} \  M_{d}(i,:)\mathbf{u} < 0 \textrm{\ or \ } u_{i} = 1 = \frac{\|\mathbf{u}\|_{1}}{\|\mathbf{u}\|^{2}_{2}}.
\end{equation}
Note that for all $i\in S$ we have $m_{ij}=1$ if and only if $j\in S$. Thus,
\begin{eqnarray*}
m_{ij}u_{j} = \left\{ \begin{array}{l} 1, \textrm{\ if\ } i, j\in S, \\
0, \textrm{\ if \ } j\notin S.
\end{array} \right.
\end{eqnarray*}
As a result, for all $i\in S$
\begin{equation}
\label{eq:Miju}
M_{d}(i,:)\mathbf{u} = \sum^{n}_{j=1} m_{ij}u_{j} = \sum_{j\in S} m_{ij}u_{j} = \sum_{j\in S}u_{j} = \sum_{j=1}^{n} u_{j} =\|\mathbf{u}\|_{1}.
\end{equation}
By combining \eqref{eq:uiM} and \eqref{eq:Miju} we obtain
\begin{equation}
\label{eq:uiMi}
u_{i} = 0 \ \textrm{and} \  M_{d}(i,:)\mathbf{u}< 0 \textrm{\quad or \quad } 1 = u_{i} = \frac{\|\mathbf{u}\|_{1}}{\|\mathbf{u}\|^{2}_{2}} = \frac{M_{d}(i,:)\mathbf{u}}{\|\mathbf{u}\|^{2}_{2}}.
\end{equation}
We can rewrite \eqref{eq:uiMi} as
\begin{equation}
\label{eq:newopticond}
\mathbf{u} = \max\left(0,\frac{M_{d}\mathbf{u}}{\|\mathbf{u}\|^{2}_{2}}\right).
\end{equation}
But \eqref{eq:newopticond} is the same as the necessary optimality condition \eqref{eq:optimalityclique2} of \eqref{eq:R1NdCG}. Therefore, $\mathbf{u} \in \textsf{S}_{\textsf{p}}$. Hence, $\mathbf{u} \in \textsf{C}_{\textsf{m}}$ if and only if $\mathbf{u}\in \textsf{F}_{\textsf{s}}\cap \textsf{S}_{\textsf{p}}$.
\qed

Theorem~\ref{thoerem4} implies that if an algorithm converges to a stationary point of~\eqref{eq:R1NdCG} and that this stationary point is binary, then it corresponds to a maximal clique.

\subsubsection{Limit Points of \eqref{eq:R1NdCG} and Feasible Solutions of \eqref{eq:MC(G)}}

This section shows how close the stationary points of \eqref{eq:R1NdCG} are to the feasible solutions of \eqref{eq:MC(G)}.

First, we present two lemmas and recall a Theorem from~\cite{GG14} about bipartite graphs. The next lemma shows that entries of $\mathbf{uu}^{\top}$ corresponding to the $-d$ entries of $M_{d}$ are approximated by zeros as $d$ gets larger.

\begin{lemma}
\label{lemmad}
For any graph $G$ and $\mathbf{u} \in \textsf{S}_{\textsf{p}}$, if $m_{ij}=-d$ and $u_{i}u_{j}>0$, we have
\begin{equation*}
0 < u_{j} < \frac{\|\mathbf{u}\|_{1}}{d+1} \quad \text{ and } \quad 0 < u_{i} < \frac{\|\mathbf{u}\|_{1}}{d+1} .
\end{equation*}
\end{lemma}
{\it Proof}
Since $u_{i}$ and $u_{j}$ are positive, the optimality condition \eqref{eq:optimalityclique2} gives
\begin{equation*}
0< u_{i}\|\mathbf{u}\|^{2}_{2} = M_{d}(i,:)\mathbf{u} = -du_{j} + \sum_{r\neq j}m_{ir}u_{r} \leq -du_{j} + \sum_{r\neq j}u_{r} = -du_{j} + \left(\|\mathbf{u}\|_{1} - u_{j}\right) = \|\mathbf{u}\|_{1} -(d+1)u_{j}.
\end{equation*}
Therefore,
$
0 < u_{j} < \frac{\|\mathbf{u}\|_{1}}{d+1}.
$
By symmetry, the same holds for $u_i$.
\qed

Below, we state and prove a lemma which is useful to draw an important relationship between stationary points of~\eqref{eq:R1NdCG} and feasible solutions of~\eqref{eq:MC(G)}.
\begin{lemma}
\label{lemma4}
Let $M$ be a symmetric matrix. If $\mathbf{u}$ is a stationary point of $\min_{\mathbf{u}\geq 0}\|M-\mathbf{uu}^\top\|^{2}_{F}$, then $(\mathbf{u},\mathbf{u})$ is a stationary point of $\min_{\mathbf{u},\mathbf{v}\geq 0}\|M-\mathbf{uv}^\top\|^{2}_{F}$.
\end{lemma}
{\it Proof}
If $\mathbf{u}=\mathbf{0}$, the proof is complete since $(\mathbf{0},\mathbf{0})$ is a stationary point of~\eqref{eq:MC(G)}.
Otherwise, since $\mathbf{u}$ is a nontrivial stationary point of $\min_{\mathbf{u}\geq 0}\|M-\mathbf{uu}^\top\|^{2}_{F}$, we have
$
\mathbf{u} = \max \left(0, \frac{M\mathbf{u}}{\|\mathbf{u}\|^{2}_{2}}\right);
$
see proof of Theorem~\ref{theorem2new}.
Moreover, if $(\mathbf{u},\mathbf{v})$
is a nontrivial stationary point of $\min_{\mathbf{u},\mathbf{v}\geq 0}\|M-\mathbf{uv}^\top\|^{2}_{F}$, by the first-order optimality conditions we have~\cite[Eq.(6)]{GG14}
\begin{equation*}
\label{eq:wh}
\mathbf{u}  = \max \left(0, \frac{M\mathbf{v}}{\|\mathbf{v}\|^{2}_{2}}\right), \qquad
\mathbf{v}  = \max \left(0, \frac{M\mathbf{u}}{\|\mathbf{u}\|^{2}_{2}}\right).
\end{equation*}
Hence, if $\mathbf{u}$ is a nontrivial stationary point of $\min_{\mathbf{u}\geq 0}\|M-\mathbf{uu}^\top\|^{2}_{F}$ then $(\mathbf{u},\mathbf{u})$ is a nontrivial stationary point of $\min_{\mathbf{u},\mathbf{v}\geq 0}\|M-\mathbf{uv}^\top\|^{2}_{F}$ (note that the converse is true only when $\mathbf{u} = \mathbf{v}$).
\qed

\begin{theorem}[Gillis and Glineur~\cite{GG14}, Th.~4 and Cor.~2]
\label{glgl}
Let $\hat{G}$ be a bipartite graph and $\hat{A}\in \{0,1\}^{m\times n}$ be its binary biadjacency matrix.
For some $\hat{d}\geq 0$, define $\hat{M}_{d} = (1+\hat{d})\hat{A} - \hat{d}\mathbf{1}_{m\times n}$, where $\mathbf{1}_{m\times n}$ is an $m$-by-$n$ matrix of ones.
Then, when $\hat{d}$ goes to infinity, every stationary point of
\begin{equation}
\label{eq:R1NdGB}
\min_{\mathbf{u}\in \mathbb{R}^{m}_{+},\mathbf{v}\in \mathbb{R}^{n}_{+}}\|\hat{M}_{d}-\mathbf{u}\mathbf{v}^{\top}\|^{2}_{F},
\end{equation}
gets arbitrarily close to some feasible solution of
\begin{equation}
\label{eq:MB}
\min_{\mathbf{u}\in \{0,1\}^{m}, \mathbf{v}\in \{0,1\}^{n}} \;
\|\hat{A} - \mathbf{u}\mathbf{v}^{\top}\|^{2}_{F}
 \quad
\text{ such that }
\quad
u_{i} + v_{j} \leq 1 + \hat{a}_{ij} \ \forall i,j.
\end{equation}
More precisely, for any $\hat{d}\geq2\max(m,n)\|\hat{A}\|_{F}$, we have for any stationary point $(\mathbf{u},\mathbf{v})$
of \eqref{eq:R1NdGB} that
\[
\min_{\mathbf{u}_{b},\mathbf{v}_{b}} \|\mathbf{u}\mathbf{v}^{\top} - \mathbf{u}_{b}\mathbf{v}_{b}^{\top}\|_{F}
< \frac{\max(m,n)\|\hat{A}\|_{F}}{\hat{d}+1},
\]
where $(\mathbf{u}_{b},\mathbf{v}_{b})$ is a feasible solution of \eqref{eq:MB}.
\end{theorem}

We can now prove our main result linking the stationary points of \eqref{eq:R1NdCG} and the feasible solutions of \eqref{eq:MC(G)}.
\begin{theorem}
\label{thmstfs}
For any graph $G$, every stationary point of \eqref{eq:R1NdCG} gets arbitrarily close to some feasible solution of \eqref{eq:MC(G)}:
\begin{equation*}
\max_{\mathbf{u}\in \textsf{S}_{\textsf{p}}}
\min_{\mathbf{u}_{c}\in \textsf{F}_{\textsf{s}}}\|\mathbf{u}-\mathbf{u}_{c}\|_{2} < \frac{n\|B\|_{F}}{d+1},
\end{equation*}
where $B$ is the modified adjacency matrix of $G$ and $d\geq 2n\|B\|_{F}$.
\end{theorem}
{\it Proof}
If the graph $G$ does not have any edge, $\textsf{S}_{\textsf{p}} = \emptyset$ and the proof is complete.
Otherwise, let $M_{d}$ be the symmetric matrix defined in \eqref{eq:matrixM}, and let $\mathbf{u} \in \textsf{S}_{\textsf{p}}$.
Then, combining Lemma~\ref{lemma4} (with $E = M_{d}$) and Theorem~\ref{glgl},
we have
\[
\min_{\mathbf{u}_{c}\in \textsf{F}_{\textsf{s}}}\|\mathbf{u}\mathbf{u}^{\top} - \mathbf{u}_{c}\mathbf{u}_{c}^{\top}\|^{2}_{F}
< \frac{n^{2}\|B\|^{2}_{F}}{(d+1)^{2}} .
\]
Observe the following:
\begin{align*}
 \min_{\mathbf{u}_{c}\in \textsf{F}_{\textsf{s}}}\|\mathbf{u}\mathbf{u}^{\top} - \mathbf{u}_{c}\mathbf{u}_{c}^{\top}\|^{2}_{F} & = \min_{\mathbf{u}_{c}\in \textsf{F}_{\textsf{s}}} \|\mathbf{u}\mathbf{u}^{\top} -\mathbf{u}_{c}\mathbf{u}^{\top} + \mathbf{u}_{c}\mathbf{u}^{\top} - \mathbf{u}_{c}\mathbf{u}_{c}^{\top}\|^{2}_{F} \\
& = \min_{\mathbf{u}_{c}\in \textsf{F}_{\textsf{s}}} \|(\mathbf{u}-\mathbf{u}_{c})\mathbf{u}^{\top} + \mathbf{u}_{c}(\mathbf{u}^{\top}-\mathbf{u}^{\top}_{c})\|^{2}_{F} \\
& \geq \min_{\mathbf{u}_{c}\in \textsf{F}_{\textsf{s}}} \|\mathbf{u}_{c}(\mathbf{u}^{\top}-\mathbf{u}^{\top}_{c})\|^{2}_{F}
= \min_{\mathbf{u}_{c}\in \textsf{F}_{\textsf{s}}} \|\mathbf{u}_{c}\|^{2}_{2}\|\mathbf{u}-\mathbf{u}_{c}\|^{2}_{2}.
\end{align*}
Therefore, since $\mathbf{u}_{c}$ is binary and $\|\mathbf{u}_{c}\|^{2}_{2}\geq 1$, we have
$
\min_{\mathbf{u}_{c}\in \textsf{F}_{\textsf{s}}}\|\mathbf{u}-\mathbf{u}_{c}\|_{2} < \frac{n\|B\|_{F}}{d+1}.
$
\qed

\begin{corollary}
\label{corr:stationarypts}
For any graph $G$, $d\geq 2n\|B\|_{F}$, and any $\mathbf{u}\in \textsf{S}_{\textsf{p}}$, we have that $\Phi(\mathbf{u}) \in \textsf{F}_{\textsf{s}}$; where
\begin{equation}
\label{eq:phiroundingop}
\Phi \ : \mathbb{R}^{n}_{+}\rightarrow \{0,1\}^{n}\ : \mathbf{u} \rightarrow \Phi(\mathbf{u}),
\end{equation}
such that for $1\leq i \leq n$
\begin{eqnarray*}
\Phi(u_{i}) = \left\{ \begin{array}{l} 0, \quad \textrm{\ if\ } u_{i}\leq 0.5, \\
1, \quad \textrm{\ if\ } u_{i}>0.5.
\end{array} \right.
\end{eqnarray*}
\end{corollary}
{\it Proof}
This follows directly from Theorem~\ref{thmstfs} since the stationary point $\mathbf{u}$ is at Euclidean distance at
most $\frac{n\|B\|_{F}}{d+1} \leq \frac{n\|B\|_{F}}{2n\|B\|_{F}+1} < \frac{1}{2}$ from a binary indicator corresponding to a clique of $G$.
\qed

\section{Proposed Algorithm and Experimental Results}
\label{sec:clqfinalg}

In this section we propose a new and efficient clique finding algorithm using our continuous formulation (Section~\ref{ouralgo}).
Our algorithm is a projected gradient scheme applied on~\eqref{eq:R1NdCG}, using the Armijo procedure for selecting the step sizes.
We then present two comparable clique finding algorithms based on the Motzkin-Straus formulation with similar computational costs (Section~\ref{sec:MSbasedalgs}).
 Finally, we provide numerical comparisons on several synthetic and real data sets in Section~\ref{exper}, where we also compare our approach to a an combinatorial heuristic proposed in~\cite{Grossoetal}.


\subsection{Projected Gradient Descent with Armijo Procedure for the Continuous Formulation~\eqref{eq:R1NdCG}} \label{ouralgo}

Consider a nonempty closed convex set $\Omega \subset \mathbb{R}^{n}$ and a continuously differentiable function $f:\mathbb{R}^{n} \rightarrow \mathbb{R}$ on $\Omega$. The projected gradient method for solving the minimization problem
\begin{equation}
\label{eq:Armminf}
\min_{\mathbf{x}\in \Omega} f(\mathbf{x}),
\end{equation}
is the following~\cite[Sec.~2.3]{Dimitri}: choose some initial $\mathbf{x}^{(1)} \in \Omega$ and, for $k = 1,2,\dots$, compute
\begin{equation*}
\label{eq:Armproj}
\mathbf{x}^{(k+1)} = \mathcal{P}_{\Omega}\left[\mathbf{x}^{(k)} - s^{(k)}\nabla f(\mathbf{x}^{(k)})\right],
\end{equation*}
where $\mathbf{x}^{(k)}$ is the $k$th iterate, $\mathcal{P}_{\Omega}$ is the projection into $\Omega$,
and $s^{(k)}$ is the step size taken at the $k$th step.
The Armijo condition requires the step $s^{(k)}$ to satisfy the condition
\begin{equation}
\label{eq:Armcond}
f(\mathbf{x}^{(k+1)}) - f(\mathbf{x}^{(k)})
\; \leq \;
\sigma \, \nabla f(\mathbf{x}^{(k)})^{\top} \left( \mathbf{x}^{(k+1)}-\mathbf{x}^{(k)} \right),
\end{equation}
for some parameter $0<\sigma < 1$.
To guarantee a sufficient decrease of the objective function at each iteration, the Armijo procedure takes $s^{(k)} = \beta^{m^{(k)}} \bar{s}$,
where $\bar{s}>0$ is a constant,
$0<\beta<1$ is a parameter,
and $m^{(k)}$ is the \emph{smallest} nonnegative integer satisfying the above condition. The limit points of a projected gradient method that uses this procedure are stationary points of~\eqref{eq:Armminf}~\cite[Prop.~2.3.3]{Dimitri}; see also~\cite{CalamaiMore}.
Searching for $s^{(k)}$ can be time consuming. Since $s^{(k-1)}$ and $s^{(k)}$ usually take values of the same order of magnitude,
using $s^{(k-1)}$ as an initial guess for $s^{(k)}$ is usually rather efficient in practice; see, e.g.,~\cite{CJLinMore}.
Algorithm~\ref{alg:clique} implements this idea on our continuous formulation~\eqref{eq:R1NdCG} of the clique problem.
The parameter $d$ in~\eqref{eq:R1NdCG} is initialized to some value
(see Section~\ref{sec:paramcond} for a discussion)
and is increased progressively (by a factor $\gamma > 1$ at each iteration)
until it reaches the upper bound $D = 2n\|A+I_{n}\|_{F} \geq n$
that guarantees
(i) the one-to-one correspondence between local and global minima of~\eqref{eq:R1NdCG} with the maximal and maximum cliques of $G$
(Theorems~\ref{theorem3} and \ref{theorem4}),
and that
(ii) rounding stationary points of~\eqref{eq:R1NdCG} gives cliques of the graph $G$ (Corollary~\ref{corr:stationarypts}).
The motivation to increase $d$ progressively is the fact that~\eqref{eq:R1NdCG} is an easy problem for small $d$.
In fact, for $d = 0$, $M_{d} = B$ is nonnegative which implies that~\eqref{eq:R1NdCG} is equivalent to computing the eigenvector of $M_{d}$ associated to the largest eigenvalue of $M_{d}$ which can be solved efficiently (combining Perron-Frobenius and Eckart-Young theorems; see, e.g.,~\cite{Golubetal}).
In fact, we have observed that trying to solve~\eqref{eq:R1NdCG} directly for a large value of $d$ leads in general to worse solutions.
\begin{algorithm}[!tb]
\caption{Clique finding algorithm}
\label{alg:clique}
\begin{algorithmic}[1]
\State $\textbf{\textrm{Require: \ }} \textrm{Adjacency matrix\ } A\in\{0,1\}^{n\times n}, \textrm{\ and\ parameters\ } \gamma>1, 0<{\sigma}<1,0<{\beta}<1$;
\State $\textbf{\textrm{Initialize: \ }} \mathbf{u}, d,\alpha$;
\State $\textbf{\textrm{Set:\ }} D = 2n\|A+I_{n}\|_{F}$;
\While{stopping criterion is not satisfied}
\State $\nabla F(\mathbf{u}) = 2 \left[\mathbf{u}\left(\|\mathbf{u}\|^{2}_{2}-1\right)-(1+d)A\mathbf{u} + d\mathbf{1}_{n}\|\mathbf{u}\|_{1}\right]$;
\State $\textrm{Fo}
= -\mathbf{u}^{\top} M_{d} \mathbf{u}
+ \frac{1}{2} \|\mathbf{u}\|^{4}_{2}
= - (1+d) \mathbf{u}^{\top} A \mathbf{u}
- (1+d) \|\mathbf{u}\|_2^2
+ d \|\mathbf{u}\|_1^2 + \frac{1}{2}\|\mathbf{u}\|^{4}_{2}$;
\While{Armijo condition is not satisfied}
\State $\mathbf{u}_{\textrm{n}} \leftarrow \max\left(0, \mathbf{u} - \alpha \nabla F(\mathbf{u})\right)$;
\State $\textrm{Fn} = -\mathbf{u}_{\textrm{n}}^{\top}M_{d}\mathbf{u}_{\textrm{n}}+ \frac{1}{2}\|\mathbf{u}_{\textrm{n}}\|^{4}_{2} = - (1+d)\mathbf{u}^{\top}_{\textrm{n}} A \mathbf{u}_{\textrm{n}} - (1+d) \|\mathbf{u}_{\textrm{n}}\|_2^2 + d \|\mathbf{u}_{\textrm{n}}\|_1^2 + \frac{1}{2}\|\mathbf{u}_{\textrm{n}}\|^{4}_{2}$;
\If{$\textrm{Fn} - \textrm{Fo} > \sigma \nabla F(\mathbf{u})^{\top} ( \mathbf{u}_{\textrm{n}}-\mathbf{u} )$}
\State $\alpha \leftarrow {\beta}\alpha$;
\Else
\State $\alpha \leftarrow \frac{\alpha}{\sqrt{{\beta}}} $;
\EndIf
\EndWhile
\State $\mathbf{u}\leftarrow \mathbf{u}_{\textrm{n}}$;
\State $d\leftarrow \min(\gamma d, D)$;
\EndWhile
\end{algorithmic}
\end{algorithm}
Note that the expressions in lines~5, ~6 and~9 of Algorithm~\ref{alg:clique} use
(i) half the objective function of \eqref{eq:matrixM} minus the constant $\|M_{d}\|_F^2$:
\[
\|M_{d} - \mathbf{u}\mathbf{u}^{\top}\|_F^2 - \|M_{d}\|_F^2 =
- 2 \mathbf{u}^{\top} M_{d} \mathbf{u} +  \|\mathbf{u}\|_2^4 ,
\]
and
(ii) the fact that the matrix $M_{d}$ is equal to $(1+d) (A + I_n) - d\mathbf{1}_{n\times n}$ so that
\begin{equation*} 
\mathbf{u}^{\top} M_{d} \mathbf{u}
=
\mathbf{u}^{\top} \left( (1+d) (A + I_n) - d\mathbf{1}_{n\times n} \right) \mathbf{u}
=
(1+d) \mathbf{u}^{\top} A \mathbf{u}
+ (1+d) \|\mathbf{u}\|_2^2
- d \|\mathbf{u}\|_1^2 .
\end{equation*}
This avoids the explicit construction of $M_{d}$, which is not practical if $A$ is sparse (since $M_{d}$ is dense).

Since every limit point of Algorithm~\ref{alg:clique} is a stationary point of~\eqref{eq:R1NdCG},
their $\Phi$ rounding, as defined in~\eqref{eq:phiroundingop}, are cliques of the graph~$G$.
\begin{theorem} \label{limptagl1}
Every limit point of Algorithm~\ref{alg:clique} is a stationary point of~\eqref{eq:R1NdCG} and the $\Phi$ rounding of these stationary points are cliques of the given graph $G$.
\end{theorem}
{\it Proof}
A projected gradient algorithm that uses the Armijo procedure converges to a stationary point~\cite{CalamaiMore}.
The second part of the theorem follows from Corollary~\ref{corr:stationarypts}, since $d$ will attain the value $D$ in a finite number of steps.
\qed

On all the numerical experiments performed in Section~\ref{exper}, Algorithm~\ref{alg:clique} always converged to a maximal clique.

\subsection{Algorithms based on the Motzkin-Straus Formulation}
\label{sec:MSbasedalgs}

In this section, we briefly describe two clique finding algorithms which are based on the Motzkin-Straus formulation \eqref{eq:MSMCP}.
The first algorithm is the relaxation scheme by Pelillo~\cite{Pelillo} and uses the following multiplicative update rule
\begin{align}
\label{eq:MSPelillo}
u^{(k+1)}_{i} = u^{(k)}_{i}\frac{(A\mathbf{u}^{(k)})_{i}}{\mathbf{u}^{(k)\top}A\mathbf{u}^{(k)}}, \tag{MSPe}
\end{align}
as an iterative procedure to find stationary points of~\eqref{eq:MSMCP},
where $A$ is the binary adjacency matrix of the graph $G$
and $\mathbf{u}^{(k)}$ is the $k$-th iterate.
The second algorithm was proposed by Ding et al.\@~\cite{Dingetal}, and generalizes the Motzkin-Straus formulation~\eqref{eq:MSMCP} as follows
\begin{equation}
\label{eq:MSDing0}
\max_{\mathbf{u}\in \mathbb{R}^{n}} \mathbf{u}^{\top}A\mathbf{u}
\textrm{\quad such that \quad }
\sum^{n}_{i=1}u_{i}^{\eta} = 1 \; \text{ and } \; \mathbf{u}\geq 0,
\end{equation}
for some $\eta \in [1,2]$ and uses the following multiplicative update rule
\begin{align}
\label{eq:MSDing1}
u^{(k+1)}_{i} = \left(u^{(k)}_{i}\frac{(A\mathbf{u}^{(k)})_{i}}{\mathbf{u}^{(k)\top}A\mathbf{u}^{(k)}}\right)^{\frac{1}{\eta}},  \tag{MSDg}
\end{align}
to find stationary points of~\eqref{eq:MSDing0}.

However, as opposed to Algorithm~\ref{alg:clique}, the above update rules do not necessarily converge to cliques of the graph $G$.
For this reason, if one wants to extract a clique from a final iterate $\mathbf{u}^*$ generated by the above updates,
a postprocessing procedure is required. The most natural strategy goes as follows:
First, the entries of the final iterate are sorted in nonincreasing order.
Then, the vertices corresponding to these sorted entries are added to a clique, one by one, until the next vertex is not connected to at least one of the vertices of the already formed clique.
For this reason, we only report the results from the post processed variants of \eqref{eq:MSPelillo} and \eqref{eq:MSDing1}.
Note that \eqref{eq:MSPelillo} was tested only on random graphs and \eqref{eq:MSDing1} was not tested on any kind of data sets by the corresponding authors. Moreover, Pelillo and Ding et al.\@ did not study the relationship between the cliques of a given graph and the stationary points of their formulation nor with the limit points of the algorithms corresponding to \eqref{eq:MSPelillo} and \eqref{eq:MSDing1}.

\subsection{Experimental Setup and Numerical Results} \label{exper}

In this section we assess the performances of the new algorithm (Algorithm~\ref{alg:clique}) compared to the two clique finding algorithms presented in section~\ref{sec:MSbasedalgs}.
 We will also compare our method with the combinatorial algorithms of Grosso et al.\@ \cite{Grossoetal} which are effective in solving the MCP. Grosso et al.\@ proposed two iterated local search algorithms based on fast neighborhood search that use multiple restarts and several thousands of node selections per second. The authors tested their algorithms on various benchmark instances.

Section~\ref{sec:paramcond} describes the parameters,  initial values and stopping criteria used for the different algorithms.
Section~\ref{sec:prepdata} describes the three types of data sets on which the experiments are run:
(i) random binary adjacency matrices,
(ii) benchmark data sets obtained from the 1992-1993 implementation challenge of DIMACS, and (iii) several text mining data sets from the CLUTO toolkit.
Finally, Section~\ref{numres} presents the numerical results, showing that Algorithm~\ref{alg:clique} outperforms the two other clique finding algorithms on these data sets.

\subsubsection{Parameters, Initial Values and Stopping Criteria}
\label{sec:paramcond}

In this section we give a brief description of the parameters, initial values and stopping criteria used in the experiments.

\paragraph{Algorithm~\ref{alg:clique}.} We use random initialization and stop the algorithm when the condition $0\leq u_{i} \leq 0.001$ or $0.999 \leq u_{i} \leq 1.001$ for all $i$ is satisfied.
In addition, we used the values ${\beta} = 0.5, {\sigma}=0.01$ and initialized $\alpha$ with $\alpha_{0} = 0.1\frac{\|\mathbf{u}^{0}\|_{2}}{\|\nabla F(\mathbf{u}^{0})\|_{2}}$ for all experiments.
To make the search for the step sizes more practical, we only try to update $\alpha$ for a maximum of five steps per iteration.
An initial value of the parameter $d$ (in Algorithm~\ref{alg:clique}) that is close to the value that balances the positive and negative entries in $M_{d}$ (that is, choosing $d$ such that $\|\max(M_{d},0)\|_F \approx \|\max(-M_{d},0)\|_F$)) works well in practice~\cite{GG14}.
For this reason, we use the initial value
$$
d_{0} = \frac{\|A+I_{n}\|^{2}_{F}}{{n^{2}-\|A+I_{n}\|^{2}_{F}}}
$$
for $d$, and then increase it by a factor $\gamma = 1.1$ at each iteration until it reaches the value $D$.
We have also experimented with SVD initialization, that is, we initialized the algorithm with the best rank-one approximation of the nonnegative modified biadjacency matrix $B=A+I_{n}$ (which is the optimal solution for $d = 0$) and found out that the results are similar to those obtained with random initialization.

\paragraph{Algorithms based on~\eqref{eq:MSPelillo} and~\eqref{eq:MSDing1}.}
 We use the initialization $\mathbf{u}^{0} = \frac{1}{n}\mathbf{1}_{n}$ and the stopping criterion $\|\mathbf{u}^{(k+1)} - \mathbf{u}^{(k)}\|^{2}_{2}<10^{-10}$ as suggested in~\cite{Pelillo}. (We also experimented with random initialization but the results were worse.)
Pelillo~\cite{Pelillo} argued that in order to avoid the algorithm \eqref{eq:MSPelillo} from being biased to a certain vertex of the graph, it is better to initialize it with values near the center of the simplex
$\Delta = \{\mathbf{u}\in \mathbb{R}^{n}_{+}: \sum^{n}_{i=1}u_{i} = 1\}$. 
For the parameter $\eta$ in \eqref{eq:MSDing1}, we used $\eta = 1.05$ as suggested in~\cite{Dingetal}. We will denote PP-\eqref{eq:MSPelillo} (resp.\@ PP-\eqref{eq:MSDing1}) the algorithm generating cliques by postprocessing the final iterate obtained with the updates~\eqref{eq:MSPelillo} (resp.\@ \eqref{eq:MSDing1}).

\paragraph{Computational cost.} Neglecting constant factors, the computational cost of these three algorithms is the same, namely $\mathcal{O}(|E|)$ per iteration, where $|E|$ is the number of edges in the graph $G$. (This applies whether the graph $G$ is sparse or not.) Most of the time is spent for computing the matrix-vector product $A\mathbf{u}$, and it can be checked that all other operations run in at most $\mathcal{O}(|V|)$ operations where $|V|$ is the number of vertices in the graph $G$.
(For example, for  Algorithm~\ref{alg:clique}, one also needs to compute the $\ell_1$ and $\ell_2$ norms of $\mathbf{u}$, $\mathbf{u}^T (A\mathbf{u})$, and to sum vectors of size $|V|$.)
For Algorithm~\ref{alg:clique}, the computational cost will be slightly higher, by a constant factor, because it requires to compute the step sizes that satisfy the Armijo condition while the two other algorithms use multiplicative updates.

\subsubsection{Preparation and Description of Data Sets}
\label{sec:prepdata}

We consider three types of data sets:
\begin{itemize}
\item[$\bullet$] Random binary matrices generated using the \texttt{sprandsym} function of \verb|MATLAB|.
For several values of density $\delta$ (the number of nonzero entries divided by the total number of entries), we generate ten $400$-by-$400$ adjacency matrices.

\item[$\bullet$] DIMACS data sets. We use the graphs available from \url{http://iridia.ulb.ac.be/~fmascia/maximum\_clique} (maintained by Franco Mascia)
where the optimal values for the clique number are available.

\item[$\bullet$] Text mining data sets from the CLUTO toolkit~\cite{Karypis}.
For convenience, we select the data sets which fit the limited memory of our personal computers.
We prepared document-by-document symmetric binary matrices as follows.
Given a document-by-word matrix (an entry $(i,j)$ is different from zero if and only if the $j$th word appears in the $i$th document),
we first set al.\@l the non-zero elements to 1, and multiply the corresponding matrix with its transpose.
Then, we set al.\@l the entries strictly less than $p-1$ to zero,
where $p$ is a parameter corresponding to the least number of words the documents are required to have in common to create an edge in the graph.
This is done for three different values of $p$: 1, 10 and 20.
Finally, we set al.\@l the nonzero entries to 1.
Binary adjacency matrices constructed in this manner correspond to graphs where the vertices are the documents, and two documents are connected to one another if they share at least $p$ words.
Hence finding cliques means extracting clusters of documents that share similar topics.

\end{itemize}

\subsubsection{Numerical Results} \label{numres}

We now present the numerical results. All tests are performed using MATLAB (R2012a) on a laptop Intel(R) Core(TM) i7-6500U CPU @2.50GHz 2.59GHz 8GB RAM. The MATLAB code is available at

\noindent \url{https://sites.google.com/site/nicolasgillis/}.

\autoref{tab:randclqsize} reports the computational cost and size of the extracted cliques for the randomly generated binary matrices. As expected, the algorithms based on \eqref{eq:MSPelillo} and~\eqref{eq:MSDing1} are in most cases computationally faster than Algorithm~\ref{alg:clique}, since they do not require to compute step sizes.
When it comes to finding cliques of larger sizes,  Algorithm~\ref{alg:clique} provides the best results in most cases.

\autoref{tab:tableDIMACS} reports the computational costs and clique sizes for the DIMACS data sets. In this experiment we included a combinatorial approach designed by Grosso et al.\@ \cite{Grossoetal} in which we directly copied the results from their paper (namely, the computational times reported in their paper are scaled to a Pentium IV 2.4 GHz with 512 MB RAM running Linux, and the value of the clique corresponds to running their algorithms 100 times and taking the average of the maximum and the minimum clique as the average clique size).
Two heuristics were proposed by Grosso et al.\@ but we considered the one with (a slightly) better performance in terms of clique size, see the results corresponding to Algorithm 2 in \textbf{Table 4} of their paper. For 13 out of 36 cases Grosso et al.\@ finds larger cliques faster than the other methods whereas on other 9 cases this method is found to be very expensive. In most cases PP-\eqref{eq:MSPelillo} and PP-\eqref{eq:MSDing1} are found to be faster than the other methods but the solutions are poor. When comparing the clique sizes of the continuous approaches, Algorithm~\ref{alg:clique} once again outperforms the other two clique finding algorithms and their post processed variants in most cases. When comparing all approaches, the discrete method due to Grosso et al.\@ provides cliques having larger sizes in many cases except in one case (hamming10\_2) where it gives the second best value (with Algorithm~\ref{alg:clique} scoring the best value) and in other five cases where there is a tie with Algorithm~\ref{alg:clique}. The success of Grosso et al.\@ can be attributed to the fact that it uses a fast neighborhood search in combination with multiple restarts which is based on selecting several thousands of nodes per second (on average 250210 nodes per second for the tested DIMACS instances). On the contrary, our algorithm is a simple single-start method which could be improved using standard techniques such as genetic algorithms or simulated annealing; for example similarly as it was done very recently for the (closely related) nonnegative matrix factorization problem~\cite{Vandaele}.  Designing such heuristics based on our approach is a direction for further research and out of the scope of this paper.

\autoref{tab:costCLUTO} reports the computational costs for extracting cliques of text mining data sets from the CLUTO toolkit. Once again, the algorithms based on \eqref{eq:MSPelillo} and~\eqref{eq:MSDing1} are computationally faster than Algorithm~\ref{alg:clique}. Note however that for some larger text data sets Algorithm~\ref{alg:clique} converges faster. However, as we will see, the extra computational cost in Algorithm~\ref{alg:clique} is worth it since it outperforms \eqref{eq:MSPelillo} and \eqref{eq:MSDing1} in almost all cases. \autoref{tab:textdatasets} contain the sizes of the cliques extracted by all the algorithms when two documents were required to have a minimum of 1/10/20 words in common to be connected in the graph, respectively.
The results show that Algorithm~\ref{alg:clique} gives the best results in all cases, except in two instances where the values of all the three algorithms coincide and in other two cases where it scores the second best value.
It is worth mentioning that Algorithm~\ref{alg:clique} converges to cliques in all the experiments, hence there is no need of postprocessing.

\begin{table}[!h]
\centering
\caption{Computational cost and clique size of graphs corresponding to $10$ randomly generated binary matrices of size $400$-by-$400$. The best results are highlighted in bold.}
\label{tab:randclqsize}
\small{
\begin{tabular}{l c l l c l l c l l}
\toprule
\ & \phantom{a} & \multicolumn{2}{l}{PP-\eqref{eq:MSPelillo}} & \phantom{a} & \multicolumn{2}{l}{PP-\eqref{eq:MSDing1}} & \phantom{a} & \multicolumn{2}{l}{Algorithm~\ref{alg:clique}} \\ \cmidrule{3-4} \cmidrule{6-7}  \cmidrule{9-10}
Density && size & time && size & time && size & time \\ \toprule
0.15 &&  \textbf{5} & {0.03} && \textbf{5} & 0.05 && \textbf{5} & \textbf{0.02} \\
0.25 &&  5 & \textbf{0.02} && \textbf{6} & 0.06 && \textbf{6} & 0.10 \\
0.35 &&  \textbf{7} & \textbf{0.03} && \textbf{7} & 0.13 && \textbf{7} & 0.11 \\
0.45 &&  7 & \textbf{0.03} && \textbf{9} & 0.22 && \textbf{9} & 0.07 \\
0.50 &&  7 & \textbf{0.05} && 9 & 0.20 && \textbf{10} & 0.13 \\
0.55 &&   8 & \textbf{0.06} && 9 & 0.23 && \textbf{10} & 0.07 \\
0.65 &&    9 & \textbf{0.03} && \textbf{11} & 0.24 && \textbf{11} & 0.09 \\
0.75 &&  11 & \textbf{0.03} && \textbf{13} & 0.25 && \textbf{13} & 0.09 \\
0.85 &&    12 & \textbf{0.03} && 13 & 0.19 && \textbf{15} & 0.09 \\
0.90 &&    13 & \textbf{0.06} && 11 & 0.20 && \textbf{14} & 0.11 \\
\bottomrule
\end{tabular}
}
\end{table}

\small{
\begin{longtable}[!h]{l l l l c l l c l l c l l c l l} 
\caption{Computational cost and clique size for selected DIMACS instances. The best and second best clique sizes are highlighted in bold and underlined, respectively.}
\label{tab:tableDIMACS}
\\ \toprule
\multicolumn{4}{l}{Test} & \phantom{a} & \multicolumn{2}{l}{PP-\eqref{eq:MSPelillo}} & \phantom{a} & \multicolumn{2}{l}{PP-\eqref{eq:MSDing1}} & \phantom{a} & \multicolumn{2}{l}{Algorithm~\ref{alg:clique}} & \phantom{a} & \multicolumn{2}{l}{Grosso et al.\@} \\ \cmidrule{1-4} \cmidrule{6-7} \cmidrule{9-10}  \cmidrule{12-13} \cmidrule{15-16}
Data set & $n$ & Edges & $\omega(G)$ && size & time && size & time && size & time && size & time\\ \toprule
\endfirsthead
\caption*{Table \ref{tab:tableDIMACS} (Continued)}\\
\toprule
\multicolumn{4}{l}{Test} & \phantom{a} & \multicolumn{2}{l}{PP-\eqref{eq:MSPelillo}} & \phantom{a} & \multicolumn{2}{l}{PP-\eqref{eq:MSDing1}} & \phantom{a} & \multicolumn{2}{l}{Algorithm~\ref{alg:clique}} & \phantom{a} & \multicolumn{2}{l}{Grosso et al.\@} \\ \cmidrule{1-4} \cmidrule{6-7} \cmidrule{9-10}  \cmidrule{12-13} \cmidrule{15-16}
Data set & $n$ & Edges & $\omega(G)$ && size & time && size & time && size & time && size & time\\ \toprule
\endhead
\bottomrule
\multicolumn{16}{r}{{Continued on next page}} \\ \bottomrule
\endfoot
\bottomrule
\endlastfoot
brock200\_1 & 200 & 14834 & 21 && 18 & 0.14 && 18 & 0.14 && \underline{19} & 0.05 && \textbf{21} & 0.02\\
brock200\_2 & 200 & 9876 & 12 && 8 & 0.03 && 9 & 0.52 && \underline{10} & 0.05 && \textbf{12} & 0.02 \\
brock200\_3 & 200 & 12048 & 15 && 10 & 0.03 && 11 & 0.09 && \underline{13} & 0.06 && \textbf{15} & 0.01 \\
brock200\_4 & 200 & 13089 & 17 && 13 & 0.03 && \underline{15} & 0.09 && \underline{15} & 0.04 && \textbf{17} & 0.13 \\
brock400\_1 & 400 & 59723 & 27 && 21 & 0.06 && 21 & 0.25 && \underline{24} & 0.08 && \textbf{27} & 9.26 \\
brock400\_2 & 400 & 59786 & 29 && 21 & 0.09 && 18 & 0.41 && \underline{24} & 0.08 && \textbf{29} & 1.20 \\
brock400\_3 & 400 & 59681 & 31 && 20 & 0.13 && 18 & 0.28 && \underline{23} & 0.09 && \textbf{31} & 0.23 \\
brock400\_4 & 400 & 59765 & 33 && 20 & 0.08 && \underline{24} & 1.03 && \underline{24} & 0.07 && \textbf{33} & 0.09 \\
brock800\_1 & 800 & 207505 & 23 && 16 & 0.19 && \underline{19} & 1.58 && 18 & 0.36 && \textbf{22.64} & 247.45 \\
brock800\_2 & 800 & 208166 & 24 && 16 & 0.14 && 16 & 1.00 && \underline{19} & 0.34 && \textbf{24} & 59.24 \\
brock800\_3 & 800 & 207333 & 25 && 18 & 0.25 && 18 & 1.77 && \underline{19} & 0.35 && \textbf{25} & 64.04 \\
brock800\_4 & 800 & 207643 & 26 && 17 & 0.19 && 16 & 1.73 && \underline{19} & 0.34 && \textbf{26} & 27.10 \\
C500-9 &  500 & 112332 & $\geq 57$ && 46 & 0.13 && 12 & 0.03 && \underline{50} & 0.12 && \textbf{57} & 1.41 \\
C1000-9 &  1000 & 450079 & $\geq 68$ && 51 & 0.75 && 5 & 0.13 && \underline{63} & 0.47 && \textbf{67.91} & 100.71 \\
C2000-5 & 2000 & 999836 & $\geq 16$ && 13 & 0.91 && \underline{14} & 8.44 && \underline{14} & 2.09 && \textbf{16} & 1.6 \\
C2000-9 &  2000 & 1799532 & $\geq 80$ && 62 & 3.86 && 7 & 0.31 && \underline{73} & 1.98 && \textbf{76.57} & 563 \\
C4000-5 & 4000 & 4000268 & $\geq 18$ && 13 & 4.47 && 3 & 1.80 && \underline{16} & 7.97 && \textbf{18} & 304.18 \\
hamming10\_2 &  1024 & 518656 & 512 && 1 & 0.001 && 1 & 0.03 && \textbf{512} & 0.13 && \underline{510.64} & 2.14 \\
keller6 & 3361 & 4619898 & $\geq 59$ && 31 & 10.17 && 15 & 1.78 && \underline{33} & 8.83 && \textbf{59} & 118.61 \\
MANN\_a27 & 378 & 70551 & 126 && 1 & 0.02 && 1 & 0.03 && \underline{123} & 0.10 && \textbf{126} & 0.005 \\
MANN\_a45 & 1035 & 533115 & 345 && 1 & 0.13 && 1 & 0.01 && \underline{333} & 0.91 && \textbf{344.02} & 373.75 \\
MANN\_a81 & 3321 & 5506380 & $\geq 1100$ && 1 & 0.95 && 1 & 0.27 && \underline{1061} & 9.24 && \textbf{1098} & 987 \\
p\_hat1000-1 &  1000 & 122253 & $\geq 10$ && 8 & 0.14 && 9 & 1.13 && \textbf{10} & 0.53 && \textbf{10} & 0.06 \\
p\_hat1000-2 &  1000 & 244799 & $\geq 46$ && 44 & 1.91 && 8 & 0.25 && \textbf{46} & 0.72 && \textbf{46} & 0.01 \\
p\_hat1000-3 &  1000 & 371746 & $\geq 68$ && \underline{63} & 0.75 && 3 & 0.13 && 62 & 0.63 && \textbf{68} & 0.07 \\
p\_hat1500-1 &  1500 & 284923 & 12 && 9 & 0.34 && \textbf{10} & 1.84 && \textbf{10} & 1.35 && \textbf{10} & 5.86 \\
p\_hat1500-2 &  1500 & 568960 & $\geq 65$ && \underline{61} & 2.09 && 22 & 0.33 && \underline{61} & 1.65 && \textbf{65} & 0.07 \\
p\_hat1500-3 &  1500 & 847244 & $\geq 94$ && 88 & 2.39 && 24 & 0.30 && \underline{92} & 1.57 && \textbf{94} & 0.09 \\
san400\_0.5\_1 & 400 & 39900 & 13 && 2 & 0.03 && 2 & 0.03 && \underline{7} & 0.04 && \textbf{13} & 0.03 \\
san400\_0.7\_1 & 400 & 55860 & 40 && 15 & 0.05 && 6 & 0.05 && \underline{22} & 0.07 && \textbf{40} & 0.04 \\
san400\_0.7\_2 & 400 & 55860 & 30 && 5 & 0.02 && 2 & 0.03 && \underline{15} & 0.04 && \textbf{30} & 0.03 \\
san400\_0.7\_3 & 400 & 55860 & 22 && 1 & 0.05 && 1 & 0.03 && \underline{13} & 0.04 && \textbf{22} & 0.05 \\
san400\_0.9\_1 & 400 & 71820 & 100 && 55 & 0.13 && 12 & 0.06 && \underline{53} & 0.04 && \textbf{100} & 0.002 \\
san1000 & 1000 & 250500 & 15 && 7 & 0.03 && 4 & 0.06 && \underline{8} & 0.23 && \textbf{15} & 2.57\\
sanr400\_0.5 & 400 & 39984 & 13 && 11 & 0.03 && 12 & 0.16 && \textbf{13} & 0.10 && \textbf{13} & 0.14\\
sanr400\_0.7 & 400 & 55869 & 21 && 18 & 0.03 && 19 & 0.25 && \textbf{21} & 0.09 && \textbf{21} & 0.02 \\
\end{longtable}
}

\begin{table}[!h]
  \centering
  \caption{Computational cost (in sec.) for finding cliques when the documents from the CLUTO toolkit were required to have at least 1/10/20 words in common, respectively. We denote $n_{\textrm{d}}$ (resp.\@ $n_{\textrm{w}}$) the number of documents (resp.\@ words) in a data set. The best results are highlighted in bold.}
\label{tab:costCLUTO}
\small{
\begin{tabular}{l l l l l l}
\toprule
Data & $n_{\textrm{d}}$ & $n_{\textrm{w}}$ & PP-\eqref{eq:MSPelillo} & PP-\eqref{eq:MSDing1} & Algorithm~\ref{alg:clique} \\
\midrule
hitech & 2,301 & 10,080 & 4.81/1.78/1.56 & \textbf{3.25/0.75/0.30} & 18.27/7.48/1.99 \\
k1b & 2,340 & 21,839 & 6.95/5.61/1.31 & \textbf{4.89/0.39/0.30} & 11.84/4.39/1.38 \\
la1 & 3,204 & 31,472 & 25.91/5.70/2.78 & 24.81/\textbf{4.27/1.27} & \textbf{11.16}/46.81/13.51 \\
la2 & 3,075 & 31,472 & 18.64/5.11/3.00 & \textbf{16.42}/\textbf{3.58/1.19} & 33.77/30.00/14.03 \\
tr23 & 204 & 5,832 & \textbf{0.02}/\textbf{0.06}/\textbf{0.03} & \textbf{0.02}/0.13/0.06 & 0.03/\textbf{0.06}/0.09 \\
tr31 & 927 & 10,127 & 0.83/0.41/0.55 & 0.75/\textbf{0.22}/\textbf{0.23} & \textbf{0.56}/1.64/0.95 \\
tr41 & 878 & 7,454 & 0.59/0.67/0.42 & \textbf{0.44/0.27/0.27} & 0.59/0.90/0.49 \\
tr45 & 690 & 8,261 & 0.38/25/0.34 & 0.27/\textbf{0.14/0.13} & \textbf{0.26}/0.90/0.47 \\
\bottomrule
\end{tabular}
}
\end{table}

\normalsize

\begin{table}[!h]
  \centering
  \caption{Clique size for the text mining data sets from the CLUTO toolkit when the documents were required to have at least 1/10/20 words in common, respectively. We denote $n_{\textrm{d}}$ (resp.\@ $n_{\textrm{w}}$) the number of documents (resp.\@ words) in a data set. The best results are highlighted in bold.}
\label{tab:textdatasets}
\small{
\begin{tabular}{l l l l l l}
\toprule
Data & $n_{\textrm{d}}$ & $n_{\textrm{w}}$ & PP-\eqref{eq:MSPelillo} & PP-\eqref{eq:MSDing1} & Algorithm~\ref{alg:clique} \\
\midrule
hitech & 2,301 & 10,080 &  1327/497/235 &  1227/458/193 & \textbf{1600/551/243} \\
k1b & 2,340 & 21,839 &  1696/302/120 &  1461/183/103 & \textbf{2040/322/122}  \\
la1 & 3,204 & 31,472 &  3077/1006/599 &  3072/956/561 & \textbf{3136/1257/639}  \\
la2 & 3,075 & 31,472 &  2627/1007/536 &  2522/961/469 & \textbf{3004/1264/612}  \\
tr23 & 204 & 5,832 &  \textbf{200}/130/\textbf{100} &  \textbf{200/134}/96 & \textbf{200}/131/99 \\
tr31 & 927 & 10,127 &  857/444/264 &  844/416/222 & \textbf{892/478/283} \\
tr41 & 878 & 7,454 &  770/391/224 &  731/351/204 & \textbf{799/412/256} \\
tr45 & 690 & 8,261 &  \textbf{689}/423/316 &  \textbf{689}/407/293 & \textbf{689/437/323} \\
\bottomrule
\end{tabular}
}
\end{table}

\normalsize

\subsection{Generalizations of Algorithm~\ref{alg:clique}}

In this paper, we focused on finding cliques in unweighted graphs.
However, Algorithm~\ref{alg:clique} can be straightforwardly used in the following two more general scenarios:
\begin{itemize}

\item[$\bullet$] \emph{Weighted graphs}. If the graph is weighted (that is, a weight is assigned to each edge of the graph indicating the importance of the relationship between two vertices), Algorithm~\ref{alg:clique} can also be used and will try to identify a clique whose corresponding submatrix has the largest possible first singular value; see formulation~\eqref{eq:MCP00}.

\item[$\bullet$] \emph{Finding dense subgraphs}. In case one is looking for dense subgraphs instead of fully connected ones, the parameter $D$ can be kept smaller. In fact, when $d$ is small, zero entries of the matrix $B$ can be approximated by positive ones. At the limit, for $d = 0$, Algorithm~\ref{alg:clique} computes the first singular vector $\mathbf{u}$ of $M_{d}$ which is positive (given that $M_{d}$ is a primitive matrix, that is, $M^{p}_{d}$ is positive for some $p$ \cite{BP94}).
The density of the graph found by Algorithm~\ref{alg:clique} will depend on the value of $D$; see also~\cite{GG14} where the idea is experimented in the case of bicliques.

\end{itemize}


\section{Conclusion}

In this paper we introduced a new continuous formulation of the maximum clique problem (MCP) using symmetric rank-one nonnegative matrix approximation; see~\eqref{eq:R1NdCG}. We showed a one-to-one correspondence between the local (resp.\@ global) optimal solutions of our continuous formulation and the maximal (resp.\@ maximum) cliques of a given graph (Theorems~\ref{theorem3} and \ref{theorem4}). In addition, we showed that the other stationary points can be made arbitrarily close to the cliques of the graph (Theorem~\ref{thmstfs}). We then proposed a new clique finding algorithm (Algorithm~\ref{alg:clique}), applying a standard projected gradient method on our continuous formulation, and showed that the limit points of this algorithm coincide with the cliques of a given graph (Theorem~\ref{limptagl1}). Finally, we tested our algorithm on various data sets: 10 randomly generated binary matrices, 36 benchmark instances from DIMACS, and 8 text mining data sets from the CLUTO toolkit. The experimental results were compared with two other continuous clique finding algorithms based on the Motzkin-Straus formulation and one discrete approach based on a fast neighborhood search that uses multiple restarts. The results show that Algorithm~\ref{alg:clique} outperforms the two continuous methods in most cases, and gives reasonable results compared to the discrete approach given that
it is a single-start local-search heuristic.


\begin{appendices}

\section{Proof of Theorem~\ref{theorem3}} \label{proofth3}

{\it Proof}
If $\mathbf{u} \in \textsf{C}_{\textsf{m}}$ then $\mathbf{u}$ is automatically binary. Then, $\mathbf{u} \in \textsf{L}_{\textsf{m}}$ if and only if there exists an $\epsilon > 0$ such that $\forall \mathbf{v} \in \mathcal{B}_{+}(\mathbf{u},\epsilon)$ we have $\left\|M_{d}-\mathbf{u}\mathbf{u}^{\top}\right\|^{2}_{F}\leq \left\|M_{d}-\mathbf{v}\mathbf{v}^{\top}\right\|^{2}_{F}$. Let $\mathbf{v}\in \mathcal{B}_{+}(\mathbf{u},\epsilon)$ and let $S_{\mathbf{u}}$ and $S_{\mathbf{v}}$ be the supports of $\mathbf{u}$ and $\mathbf{v}$, respectively. For $\epsilon < 1$, since $\mathbf{u}$ is binary, we have $S_{\mathbf{u}}\subseteq S_{\mathbf{v}}$ (that is, if $u_{i}=1$ then $v_{i}>0$).

Next, observe the following:
\begin{align}
\left\|M_{d}-\mathbf{u}\mathbf{u}^{\top}\right\|^{2}_{F} & = \left\|M_{d}(S_{\mathbf{u}},S_{\mathbf{u}})-\mathbf{u}(S_{\mathbf{u}})\mathbf{u}(S_{\mathbf{u}})^{\top}\right\|^{2}_{F} + \sum_{i,j\notin S_{\mathbf{u}}} m^{2}_{ij}\notag \\
& = \left\|M_{d}(S_{\mathbf{u}},S_{\mathbf{u}})-\mathbf{u}(S_{\mathbf{u}})\mathbf{u}(S_{\mathbf{u}})^{\top}\right\|^{2}_{F} + \sum_{i,j\in S_{\mathbf{v}}\setminus S_{\mathbf{u}}} m^{2}_{ij} + \sum_{i,j\notin S_{\mathbf{v}}} m^{2}_{ij}\notag \\
& = \left\|M_{d}(S_{\mathbf{v}},S_{\mathbf{v}})-\mathbf{u}(S_{\mathbf{v}})\mathbf{u}(S_{\mathbf{v}})^{\top}\right\|^{2}_{F} + \sum_{i,j\notin S_{\mathbf{v}}} m^{2}_{ij}, \textrm{\ and} \label{eq:muu}\\
\left\|M_{d}-\mathbf{v}\mathbf{v}^{\top}\right\|^{2}_{F}
& = \left\|M_{d}(S_{\mathbf{v}},S_{\mathbf{v}})-\mathbf{v}(S_{\mathbf{v}})\mathbf{v}(S_{\mathbf{v}})^{\top}\right\|^{2}_{F} + \sum_{i,j\notin S_{\mathbf{v}}} m^{2}_{ij}.\label{eq:mvv}
\end{align}
Therefore, by \eqref{eq:muu} and \eqref{eq:mvv}, for any $\epsilon < 1$, we have
\begin{equation*}
\left\|M_{d}-\mathbf{u}\mathbf{u}^{\top}\right\|^{2}_{F}\leq \left\|M_{d}-\mathbf{v}\mathbf{v}^{\top}\right\|^{2}_{F} \Leftrightarrow \left\|M_{d}(S_{\mathbf{v}}, S_{\mathbf{v}})-\mathbf{u}(S_{\mathbf{v}})\mathbf{u}(S_{\mathbf{v}})^{\top}\right\|^{2}_{F}\leq \left\|M_{d}(S_{\mathbf{v}},S_{\mathbf{v}})-\mathbf{v}(S_{\mathbf{v}})\mathbf{v}(S_{\mathbf{v}})^{\top}\right\|^{2}_{F}.
\end{equation*}

Let $\bar{S}_{\mathbf{u}}=S_{\mathbf{v}}\setminus S_{\mathbf{u}}$. Since $\mathbf{v}\in \mathcal{B}_{+}(\mathbf{u},\epsilon)$, there exists a $\delta \mathbf{u}$ such that $\mathbf{v} = \mathbf{u}+ \epsilon \delta \mathbf{u}$ with $\|\delta \mathbf{u}\|_{2}\leq 1$ and $\delta \mathbf{u}(\bar{S}_{\mathbf{u}})\geq 0$ since $\mathbf{u}(\bar{S}_{\mathbf{u}})=0$.

In order to express the norm $\left\|M_{d}(S_{\mathbf{v}},S_{\mathbf{v}}) - \mathbf{v}(S_{\mathbf{v}})\mathbf{v}(S_{\mathbf{v}})^{\top}\right\|^{2}_{F}$ in a more convenient way, we decompose the matrix $M_{d}(S_{\mathbf{v}}, S_{\mathbf{v}})$ into four submatrices using the decomposition $S_{\mathbf{v}} = S_{\mathbf{u}} \cup \bar{S}_{\mathbf{u}}$.
\begin{enumerate}
\item Submatrix $M_{d}(S_{\mathbf{u}}, S_{\mathbf{u}})$: Since $M_{d}(S_{\mathbf{u}}, S_{\mathbf{u}}) =\mathbf{1}_{|S_{\mathbf{u}}|\times |S_{\mathbf{u}}|}$ and $\mathbf{u}(S_{\mathbf{u}})=\mathbf{1}_{|S_{\mathbf{u}}|}$,
\begin{equation*}
e_{1} = \left\|M_{d}(S_{\mathbf{u}}, S_{\mathbf{u}}) - \mathbf{v}(S_{\mathbf{u}})\mathbf{v}(S_{\mathbf{u}})^{\top}\right\|^{2}_{F}\geq\left\|M_{d}(S_{\mathbf{u}}, S_{\mathbf{u}}) - \mathbf{u}(S_{\mathbf{u}})\mathbf{u}(S_{\mathbf{u}})^{\top}\right\|^{2}_{F}=0.
\end{equation*}
\item Submatrix $M_{d}(\bar{S}_{\mathbf{u}}, \bar{S}_{\mathbf{u}})$: Since $\mathbf{u}(\bar{S}_{\mathbf{u}})=0$ and
\begin{equation}
\label{eq:M_{d}Frobbnd}
\left\|M_{d}(\bar{S}_{\mathbf{u}}, \bar{S}_{\mathbf{u}})\right\|^{2}_{F} \leq |\bar{S}_{\mathbf{u}})|^{2}d^{2} < n^{2}(d+1)^{2}
\end{equation}
for $d\geq 1$,
\begin{align*}
e_{2} & = \left\|M_{d}(\bar{S}_{\mathbf{u}}, \bar{S}_{\mathbf{u}})-\mathbf{v}(\bar{S}_{\mathbf{u}})  \mathbf{v}(\bar{S}_{\mathbf{u}})^{\top}\right\|^{2}_{F} = \left\|M_{d}(\bar{S}_{\mathbf{u}}, \bar{S}_{\mathbf{u}})-\epsilon^{2}\delta\mathbf{u}(\bar{S}_{\mathbf{u}}) \delta\mathbf{u}(\bar{S}_{\mathbf{u}})^{\top}\right\|^{2}_{F} \\
& = \left\|M_{d}(\bar{S}_{\mathbf{u}}, \bar{S}_{\mathbf{u}})\right\|^{2}_{F} - 2\epsilon^{2}\left\langle M_{d}(\bar{S}_{\mathbf{u}}, \bar{S}_{\mathbf{u}}), \delta \mathbf{u}(\bar{S}_{\mathbf{u}}) \mathbf{u}(\bar{S}_{\mathbf{u}}) \right\rangle + \epsilon^{4}\left\|\delta \mathbf{u}(\bar{S}_{\mathbf{u}}) \delta\mathbf{u}(\bar{S}_{\mathbf{u}})\right\|^{2}_{F}\\
& \geq \left\|M_{d}(\bar{S}_{\mathbf{u}}, \bar{S}_{\mathbf{u}})\right\|^{2}_{F} - 2\epsilon^{2}\left\langle M_{d}(\bar{S}_{\mathbf{u}}, \bar{S}_{\mathbf{u}}), \delta \mathbf{u}(\bar{S}_{\mathbf{u}}) \delta\mathbf{u}(\bar{S}_{\mathbf{u}}) \right\rangle \\
& \geq \left\|M_{d}(\bar{S}_{\mathbf{u}}, \bar{S}_{\mathbf{u}})\right\|^{2}_{F} - 2\epsilon^{2}\left\|M_{d}(\bar{S}_{\mathbf{u}}, \bar{S}_{\mathbf{u}})\right\|_{F}\left\|\delta \mathbf{u}(\bar{S}_{\mathbf{u}}) \delta\mathbf{u}(\bar{S}_{\mathbf{u}})\right\|_{F} \\
& \stackrel{\eqref{eq:M_{d}Frobbnd}}{\geq} -C\epsilon ^{2} \left\|\delta\mathbf{u}(\bar{S}_{\mathbf{u}})\delta\mathbf{u}(\bar{S}_{\mathbf{u}})^{\top}\right\|_{F},
\end{align*}
where $C = 2n(d+1)$.
\item Submatrix $M_{d}(S_{\mathbf{u}}, \bar{S}_{\mathbf{u}})$: Since $\mathbf{u}(S_{\mathbf{u}})=\mathbf{1}_{|S_{\mathbf{u}}|}, \mathbf{u}(\bar{S}_{\mathbf{u}})=0, d\geq 1$ and $\epsilon < 1$,
\begin{align*}
e_{3} & = \left\|M_{d}(S_{\mathbf{u}}, \bar{S}_{\mathbf{u}})-\mathbf{v}(S_{\mathbf{u}})  \mathbf{v}(\bar{S}_{\mathbf{u}})^{\top}\right\|^{2}_{F}  = \left\|M_{d}(S_{\mathbf{u}}, \bar{S}_{\mathbf{u}})-\mathbf{v}(S_{\mathbf{u}}) \delta\mathbf{u}(\bar{S}_{\mathbf{u}})^{\top}\right\|^{2}_{F} \\
& = \left\|M_{d}(S_{\mathbf{u}}, \bar{S}_{\mathbf{u}})-\epsilon \mathbf{1}_{|S_{\mathbf{u}}|\times 1} \delta\mathbf{u}(\bar{S}_{\mathbf{u}})^{\top} - \epsilon^{2}\delta\mathbf{u}(\bar{S}_{\mathbf{u}})^{\top}\delta\mathbf{u}(S_{\mathbf{u}})^{\top}\right\|^{2}_{F} \\
& \geq \left\|M_{d}(S_{\mathbf{u}}, \bar{S}_{\mathbf{u}})-\epsilon \mathbf{1}_{|S_{\mathbf{u}}|\times 1} \delta\mathbf{u}(\bar{S}_{\mathbf{u}})^{\top}\right\|^{2}_{F} - C\epsilon^{2}\left\|\delta\mathbf{u}(\bar{S}_{\mathbf{u}})^{\top}\delta\mathbf{u}(S_{\mathbf{u}})^{\top}\right\|_{F}
\end{align*}
In fact, observe the following: Since $(\delta \mathbf{u})_{i}\leq 1, \forall i$ and $m_{ij}\in \{1,-d\}$ we have that
\begin{equation*}
\max \left(M_{d}(S_{\mathbf{u}}, \bar{S}_{\mathbf{u}})-\epsilon \mathbf{1}_{|S_{\mathbf{u}}|\times 1} \delta\mathbf{u}(\bar{S}_{\mathbf{u}})^{\top}\right)_{ij}^{2}\leq (-d-\epsilon)^{2} < (d+1)^{2}.
\end{equation*}
Therefore,
\begin{align*}
&\left\|M_{d}(S_{\mathbf{u}}, \bar{S}_{\mathbf{u}})-\epsilon \mathbf{1}_{|S_{\mathbf{u}}|\times 1} \delta\mathbf{u}(\bar{S}_{\mathbf{u}})^{\top}\right\|^{2}_{F} < |S_{u}||\bar{S}_{u}|(d+1)^{2} < n^{2}(d+1)^{2} \\
&\Rightarrow \left\|M_{d}(S_{\mathbf{u}}, \bar{S}_{\mathbf{u}})-\epsilon \mathbf{1}_{|S_{\mathbf{u}}|\times 1} \delta\mathbf{u}(\bar{S}_{\mathbf{u}})^{\top}\right\|_{F} < n(d+1).
\end{align*}
Since $\mathbf{u}$ corresponds to a maximal clique of $G$, each column of $M_{d}(S_{\mathbf{u}}, \bar{S}_{\mathbf{u}})$ must contain at least one $-d$ entry. Next we analyze each column separately. For any $i\in \bar{S}_{\mathbf{u}}$, let $n_{i}\geq 1$ be the number of $-d$ entries in the column $M_{d}(S_{\mathbf{u}},i)$ ($\therefore |S_{\mathbf{u}}|\geq n_{i}$). Then, we have
\begin{align*}
\left\|M_{d}(S_{\mathbf{u}},i) - \epsilon\mathbf{1}_{|S_{\mathbf{u}}|}(\delta \mathbf{u})_{i}\right\|^{2}_{F} & = n_{i}(-d-\epsilon (\delta \mathbf{u})_{i})^{2} + \left(|S_{\mathbf{u}}|-n_{i}\right)(1-\epsilon (\delta \mathbf{u})_{i})^{2} \\
& = n_{i}(d^{2} + 2\epsilon d (\delta \mathbf{u})_{i} + \epsilon^{2} (\delta \mathbf{u})_{i}^{2}) + (|S_{u}|-n_{i})(1-2\epsilon(\delta \mathbf{u})_{i} + \epsilon^{2} (\delta \mathbf{u})^{2}_{i}) \\
& = n_{i}d^{2} + |S_{u}|-n_{i} + 2\epsilon d \left(\delta \mathbf{u}\right)_{i} + n_{i}\epsilon^{2}
(\delta \mathbf{u})_{i}^{2}) - 2\left(|S_{u}|-n_{i}\right)\epsilon(\delta \mathbf{u})_{i} \\
& \quad + \epsilon^{2} (|S_{u}|-n_{i})(\delta \mathbf{u})^{2}_{i} \\
& \geq n_{i}d^{2} + |S_{u}|-n_{i} + 2 \epsilon (\delta \mathbf{u})_{i} \left(n_{i}(d+1) - |S_{u}|\right) \\
& \geq n_{i}d^{2} + 2 \epsilon (\delta \mathbf{u})_{i} \left(n_{i}(d+1) - |S_{u}|\right) \\
& = \left\|M_{d}(S_{\mathbf{u}},i)\right\|^{2}_{F} + 2 \epsilon (\delta \mathbf{u})_{i} \left(n_{i}(d+1) - |S_{u}|\right) \\
& \geq \|M_{d}(S_{\mathbf{u}},i)\|^{2}_{F} + 2 \epsilon (\delta \mathbf{u})_{i}.
\end{align*}
Finally, since $\delta \mathbf{u}(\bar{S}_{u})\geq 0$, summing on indices $i\in \bar{S}_{u}$ gives
\begin{align*}
e_{3}\geq \left\|M_{d}(S_{\mathbf{u}}, \bar{S}_{\mathbf{u}})-\mathbf{u}(S_{\mathbf{u}})  \mathbf{u}(\bar{S}_{\mathbf{u}})^{\top}\right\|^{2}_{F} + 2 \epsilon \left\|\delta \mathbf{u}(\bar{S}_{\mathbf{u}})\right\|_{1} - C\epsilon^{2}\left\|\delta\mathbf{u}(S_{\mathbf{u}}) \delta\mathbf{u}(\bar{S}_{\mathbf{u}})^{\top}\right\|_{F}.
\end{align*}
\item Submatrix $M_{d}(\bar{S}_{\mathbf{u}}, S_{\mathbf{u}})$: By symmetry we obtain
\begin{align*}
e_{4} & = \left\|M_{d}(\bar{S}_{\mathbf{u}}, S_{\mathbf{u}})-\mathbf{v}(\bar{S}_{\mathbf{u}})  \mathbf{v}(S_{\mathbf{u}})^{\top}\right\|^{2}_{F} \\
& \geq \left\|M_{d}(\bar{S}_{\mathbf{u}}, S_{\mathbf{u}})-\mathbf{u}(\bar{S}_{\mathbf{u}})  \mathbf{u}(S_{\mathbf{u}})^{\top}\right\|^{2}_{F} + 2 \epsilon \left\|\delta \mathbf{u}(\bar{S}_{\mathbf{u}})\right\|_{1} - C\epsilon^{2}\left\|\delta\mathbf{u}(\bar{S}_{\mathbf{u}}) \delta\mathbf{u}(S_{\mathbf{u}})^{\top}\right\|_{F}.
\end{align*}
\end{enumerate}
Combining the above results and keeping in mind that $\|\mathbf{x}\|_{2}\leq \|\mathbf{x}\|_{1}$ for any $\mathbf{x}\in \mathbb{R}^{n}$ and $\|\delta \mathbf{u}\|_{2}\leq 1$, we have that for any $0<\epsilon<\frac{1}{C}$
\begin{align*}
e_{T} & = e_{1}+ e_{2} + e_{3} + e_{4} = \left\|M_{d}(S_{\mathbf{v}}, S_{\mathbf{v}})-\mathbf{v}(S_{\mathbf{v}})  \mathbf{v}(S_{\mathbf{v}})^{\top}\right\|^{2}_{F} \\
& \geq \left\|M_{d}(S_{\mathbf{v}}, S_{\mathbf{v}})-\mathbf{u}(S_{\mathbf{v}})  \mathbf{u}(S_{\mathbf{v}})^{\top}\right\|^{2}_{F} + 4\epsilon \left\|\delta \mathbf{u}(\bar{S}_{\mathbf{u}})\right\|_{1}-2
C\epsilon^{2}\left\|\delta\mathbf{u}(\bar{S}_{\mathbf{u}})\right\|_{2} \left\|\delta\mathbf{u}(S_{\mathbf{u}})\right\|_{2} - C\epsilon^{2}\left\|\delta \mathbf{u}(\bar{S}_{\mathbf{u}})\right\|^{2}_{2} \\
& \geq \left\|M_{d}(S_{\mathbf{v}}, S_{\mathbf{v}})-\mathbf{u}(S_{\mathbf{v}})  \mathbf{u}(S_{\mathbf{v}})^{\top}\right\|^{2}_{F} + 4\epsilon \left\|\delta \mathbf{u}(\bar{S}_{\mathbf{u}})\right\|_{2}-2
C\epsilon^{2}\left\|\delta\mathbf{u}(\bar{S}_{\mathbf{u}})\right\|^{2}_{2} \left\|\delta\mathbf{u}(S_{\mathbf{u}})\right\|_{2} - C\epsilon^{2}\left\|\delta \mathbf{u}(\bar{S}_{\mathbf{u}})\right\|^{2}_{2} \\
& \geq \left\|M_{d}(S_{\mathbf{v}}, S_{\mathbf{v}})-\mathbf{u}(S_{\mathbf{v}})  \mathbf{u}(S_{\mathbf{v}})^{\top}\right\|^{2}_{F} +2\epsilon\left\|\delta \mathbf{u}(\bar{S}_{\mathbf{u}})\right\|_{2}\left(2-C\epsilon\left\|\delta \mathbf{u}(\bar{S}_{\mathbf{u}})\right\|_{2}\left\|\delta \mathbf{u}(S_{\mathbf{u}})\right\|_{2} -C\epsilon \left\|\delta \mathbf{u}(\bar{S}_{\mathbf{u}})\right\|_{2}\right) \\
& \geq \left\|M_{d}(S_{\mathbf{v}}, S_{\mathbf{v}})-\mathbf{u}(S_{\mathbf{v}})  \mathbf{u}(S_{\mathbf{v}})^{\top}\right\|^{2}_{F} +4\epsilon\left\|\delta \mathbf{u}(\bar{S}_{\mathbf{u}})\right\|_{2}(1-C\epsilon) \\
& \geq \left\|M_{d}(S_{\mathbf{v}}, S_{\mathbf{v}})-\mathbf{u}(S_{\mathbf{v}})  \mathbf{u}(S_{\mathbf{v}})^{\top}\right\|^{2}_{F}.
\end{align*}
Finally, for any $d\geq n$, $\mathbf{u}\in \textsf{C}_{\textsf{m}}$, $0<\epsilon<\frac{1}{2n(d+1)}$, $\mathbf{v}\in \mathcal{B}_{+}(\mathbf{u},\epsilon)$, we have $\left\|M_{d}-\mathbf{u}\mathbf{u}^{\top}\right\|^{2}_{F}\leq \left\|M_{d}-\mathbf{v}\mathbf{v}^{\top}\right\|^{2}_{F}$.
\qed

\end{appendices}

\section{Acknowledgment} 

We would like to thank the reviewers and Prof.\@ Kunal Narayan Chaudhury for their insightful feedback that helped us improve the paper significantly.

\end{document}